\documentstyle[epsf]{article}
\input psfig.tex
\begin{document}           
\title{On the TQFT representations of the mapping class groups
\footnote{
This preprint is available electronically at 
          \tt  http://www-fourier.ujf-grenoble.fr/\~{ }funar}\\
\small \tt (to appear Pacific Journal of Math.)}
\author{ Louis Funar\\
\small \em Institut Fourier BP 74, Univ. Grenoble I\\
\small \em 38402 Saint-Martin-d'H\`eres cedex, France\\
\small \em e-mail: funar@fourier.ujf-grenoble.fr}

\date{}
\maketitle

{\abstract We prove that the image of the mapping class group by the 
representations arising in the $SU(2)$-TQFT is infinite,
provided that the genus $g\geq 2$ and the level of the theory 
$r\neq 2,3,4,6$ (and $r\neq 10$ for $g=2$). 
In particular it follows that the quotient 
groups ${\cal M}_g/N(t^{r})$ by the normalizer of the $r$-th power
of a Dehn twist $t$ are infinite if $g\geq 3$ and $r\neq 2,3,4,6,8,12$.             }
\vspace{10mm}




\section{Introduction}
Witten \cite{Wi} constructed a TQFT in dimension 3 using path
integrals
and afterwards several rigorous constructions arose, like those using
the quantum group approach (\cite{ReTu,KiMe}), the Temperley-Lieb
algebra (\cite{L1,L2}), the theory based on the Kauffman bracket
(\cite{BHMV1,BHMV2}) or that obtained from the mapping class group
representations and the conformal field theory (\cite{Kohno}). 

Any TQFT gives rise to a tower of representations of the mapping class
groups ${\cal M}_g$ in all genera $g$ and this tower determines in
fact the theory, up to the choice of the vacuum vector (see
\cite{F1,BHMV2,Wa}). The aim of this paper is to answer whether the 
image of the mapping class groups is finite or not under such
representations. 

There is some evidence supporting the finiteness of this image group.
First, in the Abelian $U(1)$-theory the representations can be
identified with the monodromy of a system of theta functions. The
latter is explicitly computed (see e.g. \cite{F3,Go}) and  
 it is easy to see that  it factors through a finite extension (due to the
projective ambiguity) of $Sp(2g, {\bf Z}/r{\bf Z})$, where $r$ is the
level of the theory. For a general Lie group $G$, the monodromy
associated to the genus 1 surfaces may also be determined (see \cite{Jef,F3})  using some 
formulas of Kac (see \cite{Kac,KacPa}) and again it factors through a
finite extension of  $Sp(2g, {\bf Z}/r{\bf Z})$, where $r$ is now the 
shifted level. 
This has already been suggested by the fact that the 
Reshetikhin-Turaev invariant for lens spaces $L_{a,b}$ ($a \leq b$ 
running over the positive integers) takes only a finite number of
distinct values, namely for the cosets mod $r$ of $a$ and $b$ (see 
for instance \cite{LiLi}).
For low levels $r=4,6$ the whole tower of
representations was described by Blanchet and Masbaum 
in \cite{BM} and by  Wright (\cite{Wr1,Wr2}) and in
particular the images are finite groups. On the other hand all TQFTs 
are associated to conformal field theories (abb. CFT) in dimension
2 (see e.g. 
\cite{F1,Wi}) and the finiteness question appeared also in the context
of classification of rational conformal field theories. For instance
in \cite{MaSe} one asks whether the algebraic CFT have finite monodromy
(which is equivalent to our problem for some classes of 
TQFTs, as for example the $SU(2)$-TQFT). Some of the
irreducible representations of $SL(2, {\bf Z}/r{\bf Z})$ which admit
extensions as monodromies of some CFT in all genera were discussed in 
\cite{E1,E2}. Also in \cite{Ke} the action of  $SL(2, {\bf Z})$ on the 
conformal blocks was computed for all quantum doubles, and it could be
proved that the image is always finite. 
Gilmer  obtained in another way  (see \cite{Gi})  the finiteness of
the image  for $g=1$, in the $SU(2)$ theory,
result which seems to be known in the conformal field theory community, and
noticed also by M. Kontsevich.  
Meantime Stanev and Todorov \cite{ST} have a 
partial answer to this question in the case of the 4-punctured sphere, as we 
will explain below.

This is the motivation for our main result:
\newtheorem{main}{Theorem}[section]
\begin{main}
The image $\rho({\cal M}_g)$ of the mapping class group ${\cal M}_g$ 
under the representation $\rho$ arising in the $SU(2)$-TQFT
(in both the BHMV and  RT versions) 
and respectively $SO(3)$-TQFT is infinite provided that 
$g\geq 2$, $r\neq 2,3,4,6$, and if $g=2$ also $r\neq 10$. 
\end{main}
Let us mention that $\rho$ is only a projective representation and thus its 
image is well defined  up to scalar multiplication by 
 roots of unity of  order $4r$. 
To explain briefly what means the two versions BHMV and RT we recall
that the $SU(2)$-TQFT was constructed either using
the Kauffman polynomial - and this is the BHMV version from
\cite{BHMV2} - or else using the Jones polynomial - and this is the RT
version from \cite{ReTu,KiMe}. The invariants obtained for closed
3-manifolds are ``almost'' the same, but their TQFT extensions are
different.
That is the reason for considering here both of them, though as it
should be very unlikely that the mapping class group representations
do not share the same properties, in the two related  cases. 

Before we proceed let us outline the relationship with the results 
from \cite{ST}, where the Schwarz problem is considered for the 
$\widehat{su(2)}$ Knizhnik-Zamolodchikov equation. 
The authors determined  whether the image of the mapping class group of 
the 4-punctured sphere is finite, thereby 
solving a particular case of our problem, however in a slightly different 
context. It would remain to identify the following two representations of the 
mapping class groups (in arbitrary genus): 
\begin{itemize}
\item one is that arising from the conformal field theory based on the 
$\widehat{su(2)}$ Knizhnik-Zamolodchikov equation. Tsuchiya, Ueno and 
Yamada \cite{TUY} constructed the CFT using tools from algebraic geometry, 
for all Riemann surfaces.  
\item the other one is that arising in the RT-version of the 
$SU(2)$-TQFT. 
\end{itemize}

There are some naturally induced representations of braid groups in both 
approaches, which can be proved to be the same by the explicit computations
of Tsuchyia and Kanie \cite{TK}. 

Presumably the two representations of the mapping class groups  are also 
 equivalent, but a complete proof 
of this fact does not exist, on author's knowledge. 
First it should be established 
that the CFT extends to a TQFT in 3 dimensions, which is 
equivalent to control the behaviour of conformal blocks sheaves
over the compactification divisor on the moduli space of curves. 
Observe that a different and direct construction of the associated  TQFT can 
be given (\cite{Kohno,Kohno2,F5}) if we assume  the  CFT
has all the properties claimed by the physicists.
Notice that  a complete solution of that problem
would furnish an entirely geometric description of the TQFT 
following Witten's prescriptions, in which the mapping class group 
representation is the monodromy of a  projectively flat  
connection on some  vector bundle  of  non-abelian theta functions over the 
Teichmuller space. 

Thus  we cannot deduce 
directly from \cite{ST}   the finiteness of the 
mapping class group representation without assuming the previous 
unproved claim. 
Our purpose is to use instead the BHMV approach which has a simple and 
firmly established construction. 
We consider braid group representations using basically the 
monodromy of the holed spheres. The data  we obtain is similar to 
that obtained by all the other means, hence also to that from \cite{ST,TK}. 
Specifically, the idea of the proof of the main theorem is to identify a 
certain subspace of the space on which ${\cal M}_g$ acts, which is invariant
to the action of a subgroup of ${\cal M}_g$, the last being a 
 quotient of a pure braid group $P_n$, $n\geq 3$. Next we observe that
 the action of $P_n$ extends naturally to an action of the whole braid
 group $B_n$, and this it turns to factor through the Hecke algebra 
$H_n(q)$ of type $A_{n-1}$ at a root of unity $q$.
This was inspired by the computations done  by  Tsuchyia and Kanie
(\cite{TK}, see also \cite{ST}) of the monodromy in the conformal field 
theory on ${\bf  P}^1$. Now an easy modification of the  Jones 
theorem (\cite{Jones2}), 
about the generic  infiniteness of the image of $B_n$ in Hecke algebra 
representations, will settle our question.

For fixing the notations, we denote by $r$ the level, which is
supposed to be in this sequel the order of the roots of unity 
which appear in the definition of the invariants for the RT-version 
and respectively a quarter (or half) of it for the BHMV-version for 
$SU(2)$ (respectively $SO(3)$).

The groups ${\cal M}_g/N(t^r)$, quotients of ${\cal M}_g$ by 
the normalizer of a power of a Dehn twist, were previously considered
for $r=2,3$
by Humphries in \cite{Hu}, and it was shown that these are finite groups
for $r=2$ and arbitrary $g$, and infinite for $g=2$ and $r\geq 3$. 
This  solved the problem 28 asked by Birman in \cite{Bir}, p.219. 
We derive a generalization of that, to all other genera $g$, namely:
\newtheorem{main1}[main]{Corollary}
\begin{main1}
The quotient groups ${\cal M}_g/N(t^{r})$ are infinite for $g\geq 3$, 
$r\neq 2,3,4,6,8,12$. 
\end{main1}
{\em Proof:} It is well-known (see \cite{MV}, p.379) that the image of a  Dehn
twist $\rho(t)$, in some nice basis, is  a diagonal matrix
whose entries are 
$(-1)^jA^{j^2+2j}$, where $A$ is a $2r$-th root of unity. 
It follows that, for odd $r$,  $\rho(t)^{r}$, 
and for even $r$, $\rho(t)^{2r}$ respectively,  
is a scalar matrix in this particular
basis, and furthermore  it is a scalar matrix in any other  basis. 
Therefore, the image group 
$\rho({\cal M}_g)$(modulo multiplication by roots of unity of order
$4r$) is a quotient of ${\cal M}_g/N(t^{r})$ (and respectively 
${\cal M}_g/N(t^{2r})$, and now the claim follows.
 $\Box$ 

Notice that the proof given by Humphries used the Jones  representation
\cite{Jones} of ${\cal M}_2$ which arises as follows: the group 
${\cal M}_2$ is viewed as a quotient of the braid group $B_6$, and   some  
Hecke algebra representation factors throughout ${\cal M}_2$. 
For $g>2$ it is only a proper subgroup  of ${\cal M}_g$  which is a quotient of
$B_{2g+2}$, so that an extension of the Jones representation to higher 
genus is not obvious.

It seems that not only the representations have infinite image,
but the set of values  the $SU(2)$-invariant 
(at a given level $r$) takes  on the  closed 3-manifolds 
of  fixed Heegaard genus $g$, is also 
infinite. Our result does not imply this stronger statement, because
the infinite image we found comes from  a subgroup of $K\subset {\cal
  M}(F)$ of homeomorphisms of the surface extending to the handlebody. 
In fact when we  twist the gluing map of  a Heegaard splitting by an
element of $K$ we obtain a manifold  homeomorphic
to the former one. However it is very likely
that the same method could be refined to  yield this stronger statement.

We think that our theorem holds also for the case $g=2, r=10$ 
and the corollary is true more generally for $g\geq 3, r\geq 3$. 
The same ideas can be used for the $SU(N)$-TQFT to show that,
in general, the  corresponding representations of ${\cal M}_g$  have infinite images.

{\bf Acknowledgements}: This work was done during author's visit at 
Columbia University, whose hospitality is gratefully acknowledged. 
We are thankful to Roland Bacher, Joan Birman, Christian Blanchet, 
Razvan Gelca, Patrick Gilmer, 
Gregor Masbaum, Jerome Los, Vlad Sergiescu, Gretchen Wright for their
suggestions and  comments and to the referee for pointing out some 
errors in the previous version and simplifying the proof. 

\section{Preliminaries}
\subsection{Hecke algebras and Temperley-Lieb algebras}
We will outline briefly, for the sake of completeness, some basic notions
concerning the Hecke algebras (see \cite{We} for more details). 
Recall that the Hecke algebra of type $A_{n-1}$ is the algebra over
${\bf C}$ generated by $1,g_1,...,g_{n-1}$ and the following
relations:
\[ g_ig_{i+1}g_i = g_{i+1}g_ig_{i+1}, \, \, i=1,2,...,n-2, \]
\[ g_ig_j = g_jg_i, \, \, \mid i - j \mid >1, \]
\[ g_i^2 = (q-1) g_i + q, \, \, i=1,2,...,n-1,  \]
where $q\in {\bf C}-\{0\}$ is a complex parameter. Denote this algebra
by $H_n(q)$. It is known (see e.g. \cite{Bou}, p.54-55) that $H_n(q)$
is isomorphic to the group algebra ${\bf C}S_n$ of the symmetric group
$S_n$, provided that $q$ is not a root of unity. 

Notice that $H_n(q)$ is the quotient of the group algebra $CB_n$  of the braid
group $B_n$. The braid group  is usually presented as generated by 
$g_1,...,g_{n-1}$, together  with the first two relations from above. In
particular there is a natural representation of $B_n$ in $H_n(q)$.

From the quadratic relation satisfied by $g_i$ it follows that $g_i$
has at most two spectral values. For $q\neq -1$ set $f_i$ for the 
spectral projection corresponding to the eigenvalue -1; then 
$g_i = q- (1+q) f_i$, and another presentation of $H_n(q)$ can be
obtained in terms of the generators $1,f_1,...,f_{n-1}$, as follows: 
\[ f_if_{i+1}f_i - q(1+q)^{-2}f_i = f_{i+1}f_if_{i+1}-q(1+q)^{-2}f_{i+1}, \, \, i=1,2,...,n-2, \]
\[ f_if_j = f_jf_i, \, \, \mid i - j \mid >1, \]
\[ f_i^2 =  f_i, \, \, i=1,2,...,n-1.  \]

The irreducible  representations of Hecke algebras are well-known in
the case when $H_n(q)$ are semi-simple, which means that $q$ is not a 
root of unity. The structure of $H_n(q)$ at roots of unity is  
more complicated (see for instance \cite{We}) and as we will be 
concerned with  this situation precisely we introduce also some smaller 
quotients (after \cite{GHJ,Jones2}), namely  the Temperley-Lieb 
algebras. 

The algebras $A_{\beta,n}$ (following the convention
from \cite{GHJ}, section 2.8)  are generated over ${\bf C}$ by 
$1, e_1,...,e_{n-1}$ and the relations:
\[ e_ie_{i+1}e_i = e_{i}e_{i-1}e_{i}=\beta^{-1}e_i, \, \, i=1,2,...,n-2, \]
\[ e_ie_j = e_je_i, \, \, \mid i-j\mid > 1, \]
\[e_i^2=e_i, \, \, i=1,2,...,n-1. \]
Remark that $A_{\beta,n}$ is a quotient of the Hecke algebra $H_n(q)$,
for $\beta = 2+q+q^{-1}$. In fact, the image of the projector $f_i$ 
is $1-e_i$. If we replace $e_i=1-f_i$ we find a
presentation of  $A_{\beta,n}$ as the quotient of $H_n(q)$ by
adding the set of relations (see also prop.2.11.1 from \cite{GHJ},
p.123): 
\[ g_ig_{i+1}g_i+g_ig_{i+1}+g_{i+1}g_i+g_i+g_{i+1}+1=0, 
\, \mbox{ for } i=1,n-1. \]

It is known that  the algebra $A_{\beta,3}$ is semi-simple for 
 $\beta\neq 1$ (or equivalently, $r\neq 3$) and $A_{\beta,4}$ is semi-simple 
for $\beta\neq 1,2$ (equivalently $r\neq 3,4$). Moreover, when semi-simple these are 
multi-matrix algebras:  
$A_{\beta,3}= M_2({\bf C})\oplus {\bf C}$, and 
$A_{\beta,4}= M_3({\bf C})\oplus M_2({\bf C})\oplus {\bf C}$, 
 (see the theorem 2.8.5, p.98 from \cite{GHJ}). 

\subsection{Mapping class group representations}
Most of the material presented here comes from \cite{L1,R,MR}.
Let $A$ be a fixed complex number and $M$ be a compact oriented
3-manifold. The skein module $S(M)$ is the vector space generated by
the isotopy classes (rel $\partial M$) of framed links, quotiented by
the (skein) relations from figure 1. 
 
\begin{figure}
\mbox{\hspace{4.5cm}}\psfig{figure=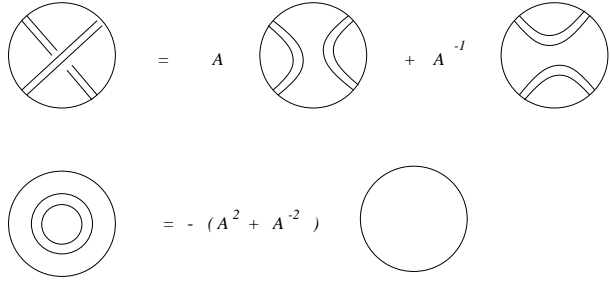,width=8cm}
\caption{Skein  relations}
\end{figure}

For example $S(S^3)$ is one dimensional (as a module over 
${\bf Z}[A,A^{-1}]$), with basis the empty link;
the image of the framed link $L\subset S^3$ in $S(S^3)$ is the value
of the Kauffman bracket evaluated at $A$. Notice that, in order to 
construct the TQFT we must specialize $A$ to be a primitive 
$2r$-th root of unity. For even  $r$   we obtain the $SU(2)$-TQFT
(level $\frac{r}{2}$ with our convention) and for odd $r$ we obtain 
the $SO(3)$-TQFT (of level $r$ this time). 

The skein space for the 3-ball with $2n$ boundary (framed) points has
an algebra structure, by representing the framed link in a planar
projection sitting into a rectangle, and separating the points into
two groups of $n$ on opposite sides. The multiplication is given by
the juxtaposition of diagrams, and the algebra  thus obtained
can be identified with the  
the Temperley-Lieb algebra  $A_{\beta,n}$, for a 
suitable $\beta$. The generators for $A_{\beta,n}$
are the elements $1_n, e_1,...,e_{n-1}$ pictured in figure 2.

\begin{figure}
\mbox{\hspace{4.5cm}}\psfig{figure=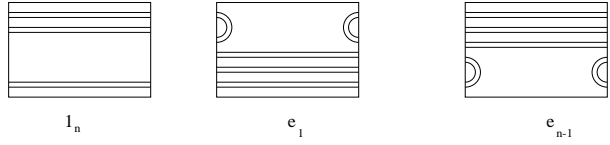,width=8cm}
\caption{Generators of the Temperley-Lieb algebra}
\end{figure}

Now the Jones-Wenzl idempotents $f^{(n)}\in A_{\beta,n}$ are uniquely  determined by
the conditions $f^{(n)2}= f^{(n)}$, $f^{(n)}e_i=e_if^{(n)}=0$, for 
$i=1,2,...,n-1$, whenever $A$ is such that all 
$\Delta_i=(-1)^i\frac{A^{2i+1}-A^{-2i-1}}{A^2-A^{-2}}$ for 
$i=0,1,...,n-1$ are non-zero. This implies that 
$f^{(n)}x=xf^{(n)}=\lambda_x f^{(n)}$, for all $x$, with a suitable 
chosen complex number $\lambda_x$.  

Denote in a planar diagram by a line labeled with $n$ (in a small
rectangle)  the element
$1_n\in A_{\beta,n}$, and by a line  with a dash labeled $n$ the insertion of
the element $f^{(n)}\in A_{\beta,n}$. This will give a convenient  description
for the elements of skein modules. 

One construction for the  $SU(2)$-invariants
via skein modules, was given in \cite{L1,L2} and latter extended to 
a TQFT in \cite{BHMV2}, and to higher $SU(N)$-invariants recently in 
\cite{L3}. 

Let us  outline first the construction of the conformal blocks, which are
the vector spaces associated to surfaces via the TQFT. Decompose the 
sphere $S^3$ as the union of two handlebodies $H$ of genus $g$, and $H'$ with a
small cylinder $F\times I$ over the surface $F=\partial H=\partial H'$
inserted between them. There is a map 
\[ <,> : S(H)\times S(H') \longrightarrow S(S^3) = {\bf C}, \]
induced by the Kauffman bracket and the union of links. In
\cite{BHMV2} it was shown that, if $A$ is a primitive $4r$-th root of
unity, then 
\[ W(F) = S(H)/\ker <,> \]
is the space associated to the surface $F$ by the  
$SU(2)$-TQFT at level $r$. Here $ker$ denotes the left
kernel of the bilinear form $<,>$. This space has however a more
concrete description. If $i,j,k$ satisfy the following 
conditions: 
\[ 0\leq i,j,k\leq r-1,\; \mid i-j\mid\leq k\leq i+j,\;  i+j+k \hbox{ is
  even}, \]  
then we can define an element of the skein space of the 3-ball with 
$i+j+k$ boundary points, given by inserting $f^{(i)}, f^{(j)},
f^{(k)}$ in the diagram, and therefore connecting up with no
crossings (see figure 3). 

\begin{figure}
\mbox{\hspace{4.5cm}}\psfig{figure=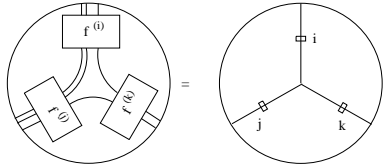,width=5cm}
\caption{Vertex elements in the skein modules}
\end{figure}

Now the triple $(i,j,k)$ is called admissible if, additionally to the
previous conditions, it satisfies $0\leq i,j,k\leq r-2$
and $i+j+k\leq 2(r-2)$. Furthermore let consider the standard 3-valent graph in
$H$ which is the standard spine of the handlebody $H$, and label its edges with
integers $i_1, i_2,..., i_{3g-3}$, such that all labels incident to a
vertex form  an admissible triple. We  form an element of $S(H)$ by
inserting idempotents $f^{(i)}$ along the edges of the graph and
triodes, like we did above at vertices.  It is shown in \cite{BHMV2,R} 
that the vectors we obtain this way form a basis of the quotient space
$W(F)$. 

For a 3-valent graph $\Gamma$, possibly with leaves and some of the
edges already carrying a label, we denote by $W(\Gamma)$ the space 
generated by the set of labelings of (non labeled) edges which have
the property that all triples from incident edges are admissible. An
easy extension of the arguments in \cite{BHMV2,R} shows that
$W(\Gamma)$ is isomorphic to $W(F)$, provided that $\Gamma$ is some 
 closed 3-valent graph of genus $g$. 

If $K$ and $K'$ are the subgroups of the mapping class group 
${\cal M}(F)$ of $F$ consisting of the classes of those homeomorphisms
which extend to the handlebodies $H$ and $H'$ respectively, then we
have natural actions of $K$ on $S(H)$, and $K'$ on $S(H')$. Moreover
these actions descend to actions on the quotient $W(F)$. One of these 
two actions, say that of $K$ on $H$, has a simple meaning: consider an
element $x\in S(H)$,  which is a representative of the class 
$[x]\in W(F)$, $u\in K$, then
 $u(x)=[\varphi(x)]\in W(F)$, where  $\varphi$ is a homeomorphism of
 $H$ whose restriction at $F$, modulo isotopy, is $u$. The other
 action, that of $K'$ on $W(F)$ can be described in a similar manner,
 using the non-degenerate bilinear form $<,>$ on $W(F)$. Namely, 
$u(x)$, for $u\in K'$, $x\in W(F)$ is defined by the equality:
\[ <ux, y> = <x, [u(y')]>, \]
holding for any $y\in W(F)$; on  the right hand side $y'\in S(H')$ is
a lift of $y$, and the action of $K'$ on $S(H')$ is the geometric one.

Moreover we have an induced action of the free group generated by $K$
and $K'$ on $W(F)$. It is shown in \cite{R,MR} that this action
descends to a central extension of the 
 mapping class group ${\cal M}(F)$. This is the
representation coming from TQFT. 
Actually we can build up the TQFT starting from that representation.
The main idea is that, if we cut a closed 3-manifold $M$ along a (closed
embedded) surface $F$ into two pieces $M_1$ and $M_2$, then 
the invariant $Z(M)$ can be recovered from the invariants $Z(M_i)$ associated
to $M_i$ (which are vectors in the space $W(F)$) as follows:
\[ Z(M) = < Z(M_1), Z(M_2)>. \]
If we want to glue back now $M_1$ to $M_2$ using an additional twist 
$\varphi \in {\cal M}(F)$ then we can compute also the invariant of
the resulting manifold $M_1\cup_{\varphi} M_2$ using the (projective)
representation $\rho: {\cal M}(F)\longrightarrow GL(W(F))$, defined
above:
\[ Z(  M_1\cup_{\varphi} M_2)= <\rho(\varphi) Z(M_1), Z(M_2)>. \]
We skipped the complications arising from the projective ambiguity, 
which is a  root of unity,  which amount to consider a supplementary 
structure  (framing) on the manifold. 
This gives a simple formula for the invariant in terms of Heegaard
splittings. In fact the vector 
$Z(H)=Z(H')\in W(F)$,
associated to the
handlebody is corresponding to the graph of genus $g$ whose labels are
all 0 (up to a normalization factor, which we skip for simplicity).
Then $Z(H\cup_{\varphi}H')$, the invariant of the closed manifold 
obtained by gluing two handlebodies along their common surface $F$
using the homeomorphism $\varphi$, is now  
$<\rho(\varphi) Z(H), Z(H')>$.

\subsection{Transformation rules for planar diagrams in the skein modules} 
In order to make explicit
computations we will freely  use  the recipes from \cite{MV} which
allows us to transform  planar diagrams representing elements in the
skein module of the 3-ball (with some boundary points) 
into simpler planar diagrams, eventually
arriving to linear combinations of the  elements of a fixed basis. 
For completeness we include these rules below.    

\vspace{0.5cm}
\begin{center} 
\parbox{2cm}{\setlength{\unitlength}{0.00050000in}%
\begingroup\makeatletter\ifx\SetFigFont\undefined%
\gdef\SetFigFont#1#2#3#4#5{%
  \reset@font\fontsize{#1}{#2pt}%
  \fontfamily{#3}\fontseries{#4}\fontshape{#5}%
  \selectfont}%
\fi\endgroup%
\begin{picture}(975,924)(826,-598)
\thicklines
\put(1651,314){\line(-2,-3){600}}
\put(1051,314){\line( 3,-5){225}}
\put(1426,-286){\line( 3,-4){225}}
\multiput(1148,-352)(10.60000,-10.60000){11}{\makebox(11.1111,16.6667){\SetFigFont{7}{8.4}{\rmdefault}{\mddefault}{\updefault}.}}
\multiput(1254,-458)(-10.60000,-10.60000){6}{\makebox(11.1111,16.6667){\SetFigFont{7}{8.4}{\rmdefault}{\mddefault}{\updefault}.}}
\multiput(1201,-511)(-10.60000,10.60000){11}{\makebox(11.1111,16.6667){\SetFigFont{7}{8.4}{\rmdefault}{\mddefault}{\updefault}.}}
\multiput(1095,-405)(10.60000,10.60000){6}{\makebox(11.1111,16.6667){\SetFigFont{7}{8.4}{\rmdefault}{\mddefault}{\updefault}.}}
\multiput(1492,-489)(10.60000,10.60000){11}{\makebox(11.1111,16.6667){\SetFigFont{7}{8.4}{\rmdefault}{\mddefault}{\updefault}.}}
\multiput(1598,-383)(10.60000,-10.60000){6}{\makebox(11.1111,16.6667){\SetFigFont{7}{8.4}{\rmdefault}{\mddefault}{\updefault}.}}
\multiput(1651,-436)(-10.60000,-10.60000){11}{\makebox(11.1111,16.6667){\SetFigFont{7}{8.4}{\rmdefault}{\mddefault}{\updefault}.}}
\multiput(1545,-542)(-10.60000,10.60000){6}{\makebox(11.1111,16.6667){\SetFigFont{7}{8.4}{\rmdefault}{\mddefault}{\updefault}.}}
\put(826, 89){\makebox(0,0)[lb]{\smash{\SetFigFont{8}{9.6}{\rmdefault}{\mddefault}{\updefault}i}}}
\put(1801, 89){\makebox(0,0)[lb]{\smash{\SetFigFont{8}{9.6}{\rmdefault}{\mddefault}{\updefault}j}}}
\end{picture}} 
\parbox{5cm}{$\; = \; \; \sum_k \delta(k;i,j)\frac{<k>}{<i,j,k>}\;   $}
\parbox{2cm}{\setlength{\unitlength}{0.00056700in}%
\begingroup\makeatletter\ifx\SetFigFont\undefined%
\gdef\SetFigFont#1#2#3#4#5{%
  \reset@font\fontsize{#1}{#2pt}%
  \fontfamily{#3}\fontseries{#4}\fontshape{#5}%
  \selectfont}%
\fi\endgroup%
\begin{picture}(825,1026)(1051,-625)
\thicklines
\multiput(1276,389)(9.00000,-9.00000){26}{\makebox(9.8039,14.7059){\SetFigFont{7}{8.4}{\rmdefault}{\mddefault}{\updefault}.}}
\multiput(1501,164)(9.00000,9.00000){26}{\makebox(9.8039,14.7059){\SetFigFont{7}{8.4}{\rmdefault}{\mddefault}{\updefault}.}}
\put(1501,164){\line( 0,-1){525}}
\multiput(1276,-586)(9.00000,9.00000){26}{\makebox(9.8039,14.7059){\SetFigFont{7}{8.4}{\rmdefault}{\mddefault}{\updefault}.}}
\multiput(1501,-361)(9.00000,-9.00000){26}{\makebox(9.8039,14.7059){\SetFigFont{7}{8.4}{\rmdefault}{\mddefault}{\updefault}.}}
\put(1651,-136){\makebox(0,0)[lb]{\smash{\SetFigFont{10}{12.0}{\rmdefault}{\mddefault}{\updefault}k}}}
\put(1051,-586){\makebox(0,0)[lb]{\smash{\SetFigFont{10}{12.0}{\rmdefault}{\mddefault}{\updefault}j}}}
\put(1876,-586){\makebox(0,0)[lb]{\smash{\SetFigFont{10}{12.0}{\rmdefault}{\mddefault}{\updefault}i}}}
\put(1126,164){\makebox(0,0)[lb]{\smash{\SetFigFont{10}{12.0}{\rmdefault}{\mddefault}{\updefault}i}}}
\put(1876,164){\makebox(0,0)[lb]{\smash{\SetFigFont{10}{12.0}{\rmdefault}{\mddefault}{\updefault}j}}}
\end{picture}}
\end{center}
\vspace{0.5cm}

\vspace{0.5cm}
\begin{center} 
\parbox{2cm}{\setlength{\unitlength}{0.00050000in}%
\begingroup\makeatletter\ifx\SetFigFont\undefined%
\gdef\SetFigFont#1#2#3#4#5{%
  \reset@font\fontsize{#1}{#2pt}%
  \fontfamily{#3}\fontseries{#4}\fontshape{#5}%
  \selectfont}%
\fi\endgroup%
\begin{picture}(675,942)(976,-766)
\thicklines
\put(1351,-286){\circle{336}}
\put(1351,164){\line( 0,-1){300}}
\put(1351,-436){\line( 0,-1){300}}
\put(976,-361){\makebox(0,0)[lb]{\smash{\SetFigFont{8}{9.6}{\rmdefault}{\mddefault}{\updefault}i}}}
\put(1426, 14){\makebox(0,0)[lb]{\smash{\SetFigFont{8}{9.6}{\rmdefault}{\mddefault}{\updefault}n}}}
\put(1201,-736){\makebox(0,0)[lb]{\smash{\SetFigFont{8}{9.6}{\rmdefault}{\mddefault}{\updefault}k}}}
\put(1651,-361){\makebox(0,0)[lb]{\smash{\SetFigFont{8}{9.6}{\rmdefault}{\mddefault}{\updefault}j}}}
\end{picture}} 
\parbox{4cm}{$\; = \; \; \delta_n^k \frac{<i,j,k>}{<k>}\;   $}
\parbox{2cm}{\setlength{\unitlength}{0.00050000in}%
\begingroup\makeatletter\ifx\SetFigFont\undefined%
\gdef\SetFigFont#1#2#3#4#5{%
  \reset@font\fontsize{#1}{#2pt}%
  \fontfamily{#3}\fontseries{#4}\fontshape{#5}%
  \selectfont}%
\fi\endgroup%
\begin{picture}(237,1149)(1489,-973)
\thicklines
\put(1576,164){\line( 0,-1){1125}}
\put(1501,-361){\line( 1, 0){150}}
\put(1651,-361){\line( 0,-1){ 75}}
\put(1651,-436){\line(-1, 0){150}}
\put(1501,-436){\line( 0, 1){ 75}}
\put(1726,-211){\makebox(0,0)[lb]{\smash{\SetFigFont{8}{9.6}{\rmdefault}{\mddefault}{\updefault}k}}}
\end{picture}
}
\end{center}
\vspace{0.5cm}

\vspace{0.5cm}
\begin{center} 
\parbox{2cm}{\setlength{\unitlength}{0.00050000in}%
\begingroup\makeatletter\ifx\SetFigFont\undefined%
\gdef\SetFigFont#1#2#3#4#5{%
  \reset@font\fontsize{#1}{#2pt}%
  \fontfamily{#3}\fontseries{#4}\fontshape{#5}%
  \selectfont}%
\fi\endgroup%
\begin{picture}(900,1149)(1351,-973)
\thicklines
\put(2101,164){\line( 0,-1){1125}}
\put(1576,164){\line( 0,-1){1125}}
\put(1501,-361){\line( 1, 0){150}}
\put(1651,-361){\line( 0,-1){ 75}}
\put(1651,-436){\line(-1, 0){150}}
\put(1501,-436){\line( 0, 1){ 75}}
\put(2026,-361){\line( 1, 0){150}}
\put(2176,-361){\line( 0,-1){ 75}}
\put(2176,-436){\line(-1, 0){150}}
\put(2026,-436){\line( 0, 1){ 75}}
\put(1351,-136){\makebox(0,0)[lb]{\smash{\SetFigFont{8}{9.6}{\rmdefault}{\mddefault}{\updefault}i}}}
\put(2251,-136){\makebox(0,0)[lb]{\smash{\SetFigFont{8}{9.6}{\rmdefault}{\mddefault}{\updefault}j}}}
\end{picture}} 
\parbox{3cm}{$\; = \; \; \sum_{k} \frac{<k>}{<i,j,k>}\;   $}
\parbox{2cm}{\setlength{\unitlength}{0.00056700in}%
\begingroup\makeatletter\ifx\SetFigFont\undefined%
\gdef\SetFigFont#1#2#3#4#5{%
  \reset@font\fontsize{#1}{#2pt}%
  \fontfamily{#3}\fontseries{#4}\fontshape{#5}%
  \selectfont}%
\fi\endgroup%
\begin{picture}(825,1026)(1051,-625)
\thicklines
\multiput(1276,389)(9.00000,-9.00000){26}{\makebox(9.8039,14.7059){\SetFigFont{7}{8.4}{\rmdefault}{\mddefault}{\updefault}.}}
\multiput(1501,164)(9.00000,9.00000){26}{\makebox(9.8039,14.7059){\SetFigFont{7}{8.4}{\rmdefault}{\mddefault}{\updefault}.}}
\put(1501,164){\line( 0,-1){525}}
\multiput(1276,-586)(9.00000,9.00000){26}{\makebox(9.8039,14.7059){\SetFigFont{7}{8.4}{\rmdefault}{\mddefault}{\updefault}.}}
\multiput(1501,-361)(9.00000,-9.00000){26}{\makebox(9.8039,14.7059){\SetFigFont{7}{8.4}{\rmdefault}{\mddefault}{\updefault}.}}
\put(1651,-136){\makebox(0,0)[lb]{\smash{\SetFigFont{10}{12.0}{\rmdefault}{\mddefault}{\updefault}k}}}
\put(1051,-586){\makebox(0,0)[lb]{\smash{\SetFigFont{10}{12.0}{\rmdefault}{\mddefault}{\updefault}j}}}
\put(1876,-586){\makebox(0,0)[lb]{\smash{\SetFigFont{10}{12.0}{\rmdefault}{\mddefault}{\updefault}i}}}
\put(1126,164){\makebox(0,0)[lb]{\smash{\SetFigFont{10}{12.0}{\rmdefault}{\mddefault}{\updefault}i}}}
\put(1876,164){\makebox(0,0)[lb]{\smash{\SetFigFont{10}{12.0}{\rmdefault}{\mddefault}{\updefault}j}}}
\end{picture}}
\end{center}
\vspace{0.5cm}

\vspace{0.5cm}
\begin{center} 
\parbox{2cm}{\setlength{\unitlength}{0.00050000in}%
\begingroup\makeatletter\ifx\SetFigFont\undefined%
\gdef\SetFigFont#1#2#3#4#5{%
  \reset@font\fontsize{#1}{#2pt}%
  \fontfamily{#3}\fontseries{#4}\fontshape{#5}%
  \selectfont}%
\fi\endgroup%
\begin{picture}(932,725)(1331,-841)
\thicklines
\put(1576,-361){\circle{474}}
\put(1801,-361){\line( 1, 0){450}}
\put(1501,-811){\makebox(0,0)[lb]{\smash{\SetFigFont{8}{9.6}{\rmdefault}{\mddefault}{\updefault}i}}}
\put(2026,-586){\makebox(0,0)[lb]{\smash{\SetFigFont{8}{9.6}{\rmdefault}{\mddefault}{\updefault}k}}}
\end{picture}
} 
\parbox{2cm}{$\; = \; \;  0   $}
\parbox{2cm}{$\;$ for $k\geq 1 \; \;    $}
\end{center}
\vspace{0.5cm}

\vspace{0.5cm}
\begin{center} 
\parbox{2cm}{\setlength{\unitlength}{0.00062500in}%
\begingroup\makeatletter\ifx\SetFigFont\undefined%
\gdef\SetFigFont#1#2#3#4#5{%
  \reset@font\fontsize{#1}{#2pt}%
  \fontfamily{#3}\fontseries{#4}\fontshape{#5}%
  \selectfont}%
\fi\endgroup%
\begin{picture}(999,825)(964,-1441)
\thicklines
\put(976,-1261){\line( 1, 1){225}}
\put(1201,-1036){\line(-1, 1){225}}
\put(1201,-1036){\line( 1, 0){525}}
\put(1951,-1261){\line(-1, 1){225}}
\put(1726,-1036){\line( 1, 1){225}}
\put(1726,-1411){\makebox(0,0)[lb]{\smash{\SetFigFont{11}{13.2}{\rmdefault}{\mddefault}{\updefault}d}}}
\put(1051,-1411){\makebox(0,0)[lb]{\smash{\SetFigFont{11}{13.2}{\rmdefault}{\mddefault}{\updefault}a}}}
\put(1726,-736){\makebox(0,0)[lb]{\smash{\SetFigFont{11}{13.2}{\rmdefault}{\mddefault}{\updefault}c}}}
\put(1126,-736){\makebox(0,0)[lb]{\smash{\SetFigFont{11}{13.2}{\rmdefault}{\mddefault}{\updefault}b}}}
\put(1426,-961){\makebox(0,0)[lb]{\smash{\SetFigFont{9}{10.8}{\rmdefault}{\mddefault}{\updefault}j}}}
\end{picture}} 
\parbox{0.5cm}{$\; = \;$}
\parbox{3cm}{$ \; \sum_{i} \; \; 
   \left\{\begin{array}{ccc}
a & b & i \\
c & d & j 
\end{array}\right\}\;\; $}
\parbox{2cm}{\setlength{\unitlength}{0.00062500in}%
\begingroup\makeatletter\ifx\SetFigFont\undefined%
\gdef\SetFigFont#1#2#3#4#5{%
  \reset@font\fontsize{#1}{#2pt}%
  \fontfamily{#3}\fontseries{#4}\fontshape{#5}%
  \selectfont}%
\fi\endgroup%
\begin{picture}(675,999)(1126,-598)
\thicklines
\put(1276,389){\line( 1,-1){225}}
\put(1501,164){\line( 1, 1){225}}
\put(1501,164){\line( 0,-1){525}}
\put(1276,-586){\line( 1, 1){225}}
\put(1501,-361){\line( 1,-1){225}}
\put(1126,-511){\makebox(0,0)[lb]{\smash{\SetFigFont{11}{13.2}{\rmdefault}{\mddefault}{\updefault}a}}}
\put(1126,164){\makebox(0,0)[lb]{\smash{\SetFigFont{11}{13.2}{\rmdefault}{\mddefault}{\updefault}b}}}
\put(1801,164){\makebox(0,0)[lb]{\smash{\SetFigFont{11}{13.2}{\rmdefault}{\mddefault}{\updefault}c}}}
\put(1276,-211){\makebox(0,0)[lb]{\smash{\SetFigFont{11}{13.2}{\rmdefault}{\mddefault}{\updefault}i}}}
\put(1726,-511){\makebox(0,0)[lb]{\smash{\SetFigFont{11}{13.2}{\rmdefault}{\mddefault}{\updefault}d}}}
\end{picture}}
\end{center}
\vspace{0.5cm}

\vspace{0.5cm}
\begin{center} 
\parbox{2cm}{\setlength{\unitlength}{0.00050000in}%
\begingroup\makeatletter\ifx\SetFigFont\undefined%
\gdef\SetFigFont#1#2#3#4#5{%
  \reset@font\fontsize{#1}{#2pt}%
  \fontfamily{#3}\fontseries{#4}\fontshape{#5}%
  \selectfont}%
\fi\endgroup%
\begin{picture}(1149,1092)(1189,-691)
\thicklines
\put(1726,389){\line( 0,-1){375}}
\put(1726, 14){\line(-3,-4){225}}
\put(1501,-286){\line( 0, 1){  0}}
\put(1501,-286){\line( 1, 0){450}}
\put(1951,-286){\line(-3, 4){225}}
\put(1501,-286){\line(-4,-3){300}}
\put(1951,-286){\line( 5,-3){375}}
\put(1801,239){\makebox(0,0)[lb]{\smash{\SetFigFont{8}{9.6}{\rmdefault}{\mddefault}{\updefault}A}}}
\put(1426,-136){\makebox(0,0)[lb]{\smash{\SetFigFont{8}{9.6}{\rmdefault}{\mddefault}{\updefault}B}}}
\put(1951,-136){\makebox(0,0)[lb]{\smash{\SetFigFont{8}{9.6}{\rmdefault}{\mddefault}{\updefault}E}}}
\put(1651,-511){\makebox(0,0)[lb]{\smash{\SetFigFont{8}{9.6}{\rmdefault}{\mddefault}{\updefault}D}}}
\put(2026,-661){\makebox(0,0)[lb]{\smash{\SetFigFont{8}{9.6}{\rmdefault}{\mddefault}{\updefault}C}}}
\put(1351,-661){\makebox(0,0)[lb]{\smash{\SetFigFont{8}{9.6}{\rmdefault}{\mddefault}{\updefault}F}}}
\end{picture}} 
\parbox{0.5cm}{$\; = \;$}
\parbox{4cm}{$\;  \frac{\left < 
                \begin{array}{ccc}
                 A & B & E \\
                 D& C  & F
                 \end{array} \right >}{<A,F,C>}\;   $}
\parbox{2cm}{\setlength{\unitlength}{0.00050000in}%
\begingroup\makeatletter\ifx\SetFigFont\undefined%
\gdef\SetFigFont#1#2#3#4#5{%
  \reset@font\fontsize{#1}{#2pt}%
  \fontfamily{#3}\fontseries{#4}\fontshape{#5}%
  \selectfont}%
\fi\endgroup%
\begin{picture}(987,924)(1051,-523)
\thicklines
\put(1576,389){\line( 0,-1){525}}
\put(1576,-136){\line(-1,-1){375}}
\put(1576,-136){\line( 6,-5){450}}
\put(1726, 89){\makebox(0,0)[lb]{\smash{\SetFigFont{8}{9.6}{\rmdefault}{\mddefault}{\updefault}A}}}
\put(1051,-436){\makebox(0,0)[lb]{\smash{\SetFigFont{8}{9.6}{\rmdefault}{\mddefault}{\updefault}F}}}
\put(2026,-436){\makebox(0,0)[lb]{\smash{\SetFigFont{8}{9.6}{\rmdefault}{\mddefault}{\updefault}C}}}
\end{picture}}
\end{center}
\vspace{0.5cm}

\vspace{0.5cm}
\begin{center} 
\parbox{2cm}{\setlength{\unitlength}{0.00062500in}%
\begingroup\makeatletter\ifx\SetFigFont\undefined%
\gdef\SetFigFont#1#2#3#4#5{%
  \reset@font\fontsize{#1}{#2pt}%
  \fontfamily{#3}\fontseries{#4}\fontshape{#5}%
  \selectfont}%
\fi\endgroup%
\begin{picture}(999,1080)(2839,-757)
\thicklines
\put(2851,-211){\line( 1, 0){975}}
\multiput(3826,-211)(-6.25000,-9.37500){25}{\makebox(8.8889,13.3333){\SetFigFont{7}{8.4}{\rmdefault}{\mddefault}{\updefault}.}}
\multiput(2851,-211)(6.25000,-9.37500){25}{\makebox(8.8889,13.3333){\SetFigFont{7}{8.4}{\rmdefault}{\mddefault}{\updefault}.}}
\multiput(3676,-436)(-10.71429,-3.57143){22}{\makebox(8.8889,13.3333){\SetFigFont{7}{8.4}{\rmdefault}{\mddefault}{\updefault}.}}
\multiput(2851,-211)(6.25000,9.37500){25}{\makebox(8.8889,13.3333){\SetFigFont{7}{8.4}{\rmdefault}{\mddefault}{\updefault}.}}
\multiput(3826,-211)(-6.25000,9.37500){25}{\makebox(8.8889,13.3333){\SetFigFont{7}{8.4}{\rmdefault}{\mddefault}{\updefault}.}}
\multiput(3001, 14)(10.71429,3.57143){22}{\makebox(8.8889,13.3333){\SetFigFont{7}{8.4}{\rmdefault}{\mddefault}{\updefault}.}}
\multiput(3676, 14)(-10.71429,3.57143){22}{\makebox(8.8889,13.3333){\SetFigFont{7}{8.4}{\rmdefault}{\mddefault}{\updefault}.}}
\put(3226, 89){\line( 1, 0){225}}
\multiput(3001,-436)(10.71429,-3.57143){22}{\makebox(8.8889,13.3333){\SetFigFont{7}{8.4}{\rmdefault}{\mddefault}{\updefault}.}}
\put(3226,-511){\line( 1, 0){225}}
\put(3301,239){\makebox(0,0)[lb]{\smash{\SetFigFont{11}{13.2}{\rmdefault}{\mddefault}{\updefault}a}}}
\put(3301,-136){\makebox(0,0)[lb]{\smash{\SetFigFont{11}{13.2}{\rmdefault}{\mddefault}{\updefault}b}}}
\put(3301,-736){\makebox(0,0)[lb]{\smash{\SetFigFont{11}{13.2}{\rmdefault}{\mddefault}{\updefault}c}}}
\end{picture}} 
\parbox{4cm}{$\; = \; \; <a,b,c>\;   $}
\end{center}
\vspace{0.5cm}

\vspace{0.5cm}
\begin{center} 
\parbox{4cm}{\setlength{\unitlength}{0.00050000in}%
\begingroup\makeatletter\ifx\SetFigFont\undefined%
\gdef\SetFigFont#1#2#3#4#5{%
  \reset@font\fontsize{#1}{#2pt}%
  \fontfamily{#3}\fontseries{#4}\fontshape{#5}%
  \selectfont}%
\fi\endgroup%
\begin{picture}(1824,1992)(364,-1891)
\thicklines
\put(1276, 89){\line( 0,-1){1200}}
\put(1276, 89){\line(-3,-2){900}}
\put(376,-511){\line( 3,-2){900}}
\put(1276,-1111){\line( 4, 3){900}}
\put(2176,-436){\line(-5, 3){893.382}}
\put(376,-511){\line( 1,-5){150}}
\put(526,-1261){\line( 3,-1){742.500}}
\put(1276,-1486){\line( 3, 1){675}}
\put(1951,-1261){\line( 1, 4){207.353}}
\put(1276,-1861){\makebox(0,0)[lb]{\smash{\SetFigFont{8}{9.6}{\rmdefault}{\mddefault}{\updefault}F}}}
\put(601,-136){\makebox(0,0)[lb]{\smash{\SetFigFont{8}{9.6}{\rmdefault}{\mddefault}{\updefault}B}}}
\put(901,-736){\makebox(0,0)[lb]{\smash{\SetFigFont{8}{9.6}{\rmdefault}{\mddefault}{\updefault}D}}}
\put(1876,-61){\makebox(0,0)[lb]{\smash{\SetFigFont{8}{9.6}{\rmdefault}{\mddefault}{\updefault}A}}}
\put(1426,-586){\makebox(0,0)[lb]{\smash{\SetFigFont{8}{9.6}{\rmdefault}{\mddefault}{\updefault}E}}}
\put(1651,-1111){\makebox(0,0)[lb]{\smash{\SetFigFont{8}{9.6}{\rmdefault}{\mddefault}{\updefault}C}}}
\end{picture}} 
\parbox{0.5cm}{$\; = \;$}
\parbox{5cm}{$\;  \left < 
                \begin{array}{ccc}
                 A & B & E \\
                 D& C  & F
                 \end{array} \right >\;   $}
\end{center}
\vspace{0.5cm}

\noindent where
\[ <k> = (-1)^{k}[k+1] = (-1)^k\frac{A^{2k+2}-A^{-2k-2}}{A^2-A^{-2}},
\]
\[ \delta(c;a,b)= (-1)^kA^{ij-k(i+j+k+2)}, \]
\[ <a,b,c> = (-1)^{i+j+k}\frac{[i+j+k+1]![i]![j]![k]!}{[i+j]![i+k]![j+k]!}.\]
Here $i$, $j$, $k$ are the internal colors given by 
\[ i=\frac{b+c-a}{2},\;  j=\frac{a+c-b}{2},\; k=\frac{b+a-c}{2}, \]
and $[n]!=[1][2]...[n]$. 

Consider now $A,B,C,D,E,F$ such that $(A,B,E)$, $(B,D,F)$, $(E,D,C)$, 
$(A,C,F)$ are admissible triples and  make some notations:
$\Sigma= A+B+C+D+E+F$, $a_1=\frac{A+B+E}{2}$,  $a_2=\frac{B+D+F}{2}$, 
$a_3=\frac{E+D+C}{2}$,  $a_4=\frac{A+C+F}{2}$, $b_1 =\frac{\Sigma-A-F}{2}$,
$b_2=\frac{\Sigma-B-C}{2}$,  $a_1=\frac{\Sigma-A-D}{2}$.

The tetrahedron coefficient is defined as:
\[\left < 
                \begin{array}{ccc}
                 A & B & E \\
                 D& C  & F
                 \end{array} \right >
=
\frac{\prod_i\prod_j [b_i-a_j]!}{[A]![B]![C]![D]![E]![F]!}
\left ( 
                \begin{array}{ccccccc}
                 a_1 &    & a_2 &     & a_3 &    & a_4 \\
                     & b_1&     & b_2 &     & b_3&
                 \end{array} \right ), \]
where
\[ \left ( 
                \begin{array}{ccccccc}
                 a_1 &    & a_2 &     & a_3 &    & a_4 \\
                     & b_1&     & b_2 &     & b_3&
                 \end{array} \right )
=
\sum_{\max{a_j}\leq \zeta\leq\min{a_j}}
\frac{(-1)^{\zeta}[\zeta +1]!}{\prod_i[b_i-\zeta]!\prod_j[\zeta-a_i]!}.
\]
The quantum 6j-symbol of \cite{MV} is given by the formula:
\[\left\{\begin{array}{ccc}
a & b & i \\
c & d & j 
\end{array}\right\}
=
\frac{<i>
      \left < 
                \begin{array}{ccc}
                 i & b & c \\
                 j& d  & a
                 \end{array} \right >}
      {<i,a,d><i,b,c>}.
\]

\section{Proof of the theorem} 
\subsection{Outline}
Consider a surface $F$ of genus $g$ and let $\Gamma \subset H$ be a 3-valent
graph embedded in the handlebody $H$. Suppose that the graph $\Gamma'$
shown in figure 4 is a subgraph of $\Gamma$. Then $\Gamma'$ can be
viewed as the spine of the $(n+2)$-holed sphere $F'\subset F$, which
is the intersection of a regular neighborhood of $\Gamma'$ (in 
${\bf R}^3$) with $F$. 
Consider a partial labeling $\Gamma'(n,m)$
of $\Gamma'$ as shown on the right hand side in the figure. Notice that the leaf with label 0
can be removed without affecting the space $W(\Gamma'(n,m))$. 
\newtheorem{lem}{Lemma}[section]
\begin{lem}
For a suitably chosen $\Gamma$, of genus $g\geq 4$, there exist 
$m\geq 0, n\geq 5$, such that $W(\Gamma'(n,m))\subset W(\Gamma)$ and 
$\dim W(\Gamma'(n,m))\geq 2$. For $g=3$ we have 
$W(\Gamma'(4,2))\subset W(\Gamma)$, and for $g=2$
$W(\Gamma'(3,1))\subset W(\Gamma)$. 
\end{lem}
{\em Proof:} We connect among them  the leaves of $\Gamma'$ 
 using some planar arcs, in order to obtain a closed graph of 
minimal genus, and such that the labels agree when making connections. $\Box$

\begin{figure}
\mbox{\hspace{4.5cm}}\psfig{figure=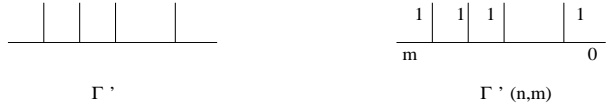,width=8cm}
\caption{The graphs $\Gamma'(n,m)$}
\end{figure}

Fix now once for all the embedding of graphs $\Gamma'\subset \Gamma$
as in the lemma, and denote by $V=V(n,m)\subset W(\Gamma)$ the image of 
$W(\Gamma'(n,m))$. Consider the curves $\gamma_{ij}\subset F'$, $1\leq
i,j\leq n$, drawn on
the $(n+2)$-holed sphere $F'$ which encircle the holes $i$ and $j$
like in the picture 5 and the set of curves $\gamma_i$ which are loops
around the holes. The Dehn twists $T_{\gamma_{ij}}$ and $T_{\gamma_i}$
generate a subgroup $S$ of the mapping class group ${\cal M}(F)$. 

\begin{figure}
\mbox{\hspace{4.5cm}}\psfig{figure=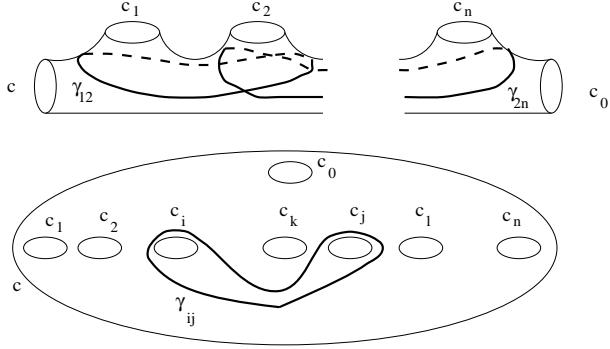,width=8cm}
\caption{The curves $\gamma_{ij}$}
\end{figure}

\newtheorem{prop}{Proposition}[section]
\begin{prop}\label{pr}
The subspace $V(n,m)\subset W(F)$ is $\rho(S)$-invariant. Moreover the
image $\rho\mid_{V(3,1)}(S)\subset GL(V(3,1))$ is an infinite group, provided
that $g\geq 2$ and the level $r\neq 2,3,4,6,10$.
For $r=10$ and $g\geq 3$  the image of 
$\rho\mid_{V(4,2)}(S)\subset GL(V(4,2))$ is infinite. 
\end{prop}
The first part of the proposition follows from  a more
general fact concerning sub-surfaces $F'\subset F$ and subgroups $S$ of 
${\cal M}(F)$ of classes of homeomorphisms  which keep $F'$ invariant
up to an isotopy, and send each boundary component into itself. Assume 
that a labeling of the boundary components of $F'$ was fixed: 
this amounts to fix a labeling for the leaves of the subgraph 
$\Gamma'$, the spine of $F'$. Then the subspace $W(\Gamma')\subset
W(\Gamma)$ is invariant by the action of $S$ on $W(F)$. Moreover,
consider now that $F'$ may be sent by a larger group $S'$ into a
sub-surface $F''$ which is isotopic to $F$, but the boundary components
may be permuted among themselves. We claim now, that a 
subspace $W(\Gamma')$, associated to a labeling of the boundary
components, is sent by such a homeomorphism into the subspace
$W(\Gamma')$ associated to the permuted labeling on the boundary
components. In particular the space $V$ from the proposition is not only
invariant under $\rho(S)$, but also under larger groups which could permute 
the $n$ boundary components $c_i$, $i=1,2,...,n$, since all their
labels are identical. 

Another observation is that the action of ${\cal M}(F')$ on the space 
$W(\Gamma')$, where $\Gamma'$ has one external edge $e$ (corresponding
to the boundary component $c_e\subset \partial F'$) is the same as the
action of  ${\cal M}(F'\cup_{c_e}D^2)$ on the space $W(\Gamma'')$;
here $F'\cup_{c_e}D^2$ is the result of gluing a disk on the circle
$c_e$, and $\Gamma''$ is $\Gamma$ with the edge $e$ removed from it.
This way we see that  $\rho(S)$ acts like 
${\cal M}(F'\cup_{c_0}D^2)$ on the given space. This will help to
find out the corresponding extension to the braid group. 

Before we  explain this action, remark that all Dehn
twists along $\gamma_{ij}, \gamma_i$ are elements of the subgroup $K\subset
{\cal M}(F)$ of classes of homeomorphisms extending to the handlebody
$H$. Therefore, according to the discussion in the previous section, 
the action of $T_{\gamma_{ij}}$ (or $T_{\gamma_i}$) on $V(n,m)$ has a simple expression in the
skein module of the 3-ball with $(n+2)$-boundary points: just perform
the Dehn twist on the 3-ball which is a regular neighborhood of the 
graph  $\Gamma'$, viewed as part of the handlebody $H$, whose spine is
$\Gamma$. This is equivalent with twisting the $i$-th and $j$-th legs of
the graph $\Gamma'$. We have to apply further the skein relations, in
order to compute the latter  element in terms of the given  basis of
$V$, where all the legs are straight. 
Notice that the representative  graphs considered 
in the sections 2.2 and 2.3 are framed graphs, and the framings 
considered in the planar pictures are the blackboard ones. When the
Dehn twist $T_{\gamma_{ij}}$ acts on the spine graph $\Gamma'$ 
then the framings of the strands $i$ and $j$ are altered. Then 
the action of $T_{\gamma_i}$ on the $i$-th strand is the  
change of its framing by one unit. Then the element 
$A_{ij}=T_{\gamma_{ij}}T_{\gamma_i}^{-1}T_{\gamma_j}^{-1}$ acts on 
$V(n,m)$ as follows: 

\vspace{0.5cm}
\begin{center}
\setlength{\unitlength}{0.00033300in}%
\begingroup\makeatletter\ifx\SetFigFont\undefined%
\gdef\SetFigFont#1#2#3#4#5{%
  \reset@font\fontsize{#1}{#2pt}%
  \fontfamily{#3}\fontseries{#4}\fontshape{#5}%
  \selectfont}%
\fi\endgroup%
\begin{picture}(11799,1620)(439,-1069)
\thicklines
\put(6976,-436){\line( 1, 0){5250}}
\put(12226,-436){\line( 0, 1){  0}}
\put(10651, 89){\line( 0,-1){525}}
\put(451,-436){\line( 1, 0){5250}}
\put(5701,-436){\line( 0, 1){  0}}
\put(1201,-436){\line( 0, 1){975}}
\put(1651,-436){\line( 0, 1){975}}
\put(2101,-436){\line( 0, 1){975}}
\put(5176,-436){\line( 0, 1){975}}
\put(4801,-436){\line( 0, 1){975}}
\put(8176,-436){\line( 0, 1){975}}
\put(9451,-436){\line( 0, 1){975}}
\put(11701,-436){\line( 0, 1){975}}
\put(11026,-436){\line( 0, 1){975}}
\put(10651,239){\line( 0, 1){300}}
\put(5851,-136){\vector( 1, 0){900}}
\put(8551,-436){\line( 0, 1){150}}
\multiput(8551,-286)(15.00000,15.00000){6}{\makebox(16.6667,25.0000){\SetFigFont{7}{8.4}{\rmdefault}{\mddefault}{\updefault}.}}
\put(8626,-211){\line( 1, 0){750}}
\put(9526,-211){\line( 1, 0){600}}
\put(10126,-211){\line( 1, 0){225}}
\multiput(10351,-211)(20.45455,6.81818){12}{\makebox(16.6667,25.0000){\SetFigFont{7}{8.4}{\rmdefault}{\mddefault}{\updefault}.}}
\put(9526,164){\line( 1, 0){825}}
\put(10351,164){\line( 1, 0){300}}
\multiput(10651,164)(18.75000,-9.37500){9}{\makebox(16.6667,25.0000){\SetFigFont{7}{8.4}{\rmdefault}{\mddefault}{\updefault}.}}
\put(10801, 89){\line( 0,-1){ 75}}
\multiput(10801, 14)(-15.00000,-15.00000){6}{\makebox(16.6667,25.0000){\SetFigFont{7}{8.4}{\rmdefault}{\mddefault}{\updefault}.}}
\put(8551,539){\line( 0,-1){225}}
\multiput(8551,314)(15.00000,-15.00000){11}{\makebox(16.6667,25.0000){\SetFigFont{7}{8.4}{\rmdefault}{\mddefault}{\updefault}.}}
\put(8701,164){\line( 1, 0){675}}
\put(7651,-436){\line( 0, 1){975}}
\put(11851,-886){\makebox(0,0)[lb]{\smash{\SetFigFont{8}{9.6}{\rmdefault}{\mddefault}{\updefault}0}}}
\put(7126,-886){\makebox(0,0)[lb]{\smash{\SetFigFont{8}{9.6}{\rmdefault}{\mddefault}{\updefault}m}}}
\put(5326,-886){\makebox(0,0)[lb]{\smash{\SetFigFont{8}{9.6}{\rmdefault}{\mddefault}{\updefault}0}}}
\put(601,-886){\makebox(0,0)[lb]{\smash{\SetFigFont{8}{9.6}{\rmdefault}{\mddefault}{\updefault}m}}}
\put(1276,-886){\makebox(0,0)[lb]{\smash{\SetFigFont{8}{9.6}{\rmdefault}{\mddefault}{\updefault}p}}}
\put(1351,-1036){\makebox(0,0)[lb]{\smash{\SetFigFont{6}{7.2}{\rmdefault}{\mddefault}{\updefault}n-1}}}
\put(1051,164){\makebox(0,0)[lb]{\smash{\SetFigFont{8}{9.6}{\rmdefault}{\mddefault}{\updefault}1}}}
\put(1501,164){\makebox(0,0)[lb]{\smash{\SetFigFont{8}{9.6}{\rmdefault}{\mddefault}{\updefault}1}}}
\put(1951,164){\makebox(0,0)[lb]{\smash{\SetFigFont{8}{9.6}{\rmdefault}{\mddefault}{\updefault}1}}}
\put(5251, 89){\makebox(0,0)[lb]{\smash{\SetFigFont{8}{9.6}{\rmdefault}{\mddefault}{\updefault}1}}}
\put(4876, 89){\makebox(0,0)[lb]{\smash{\SetFigFont{8}{9.6}{\rmdefault}{\mddefault}{\updefault}1}}}
\put(5026,-961){\makebox(0,0)[lb]{\smash{\SetFigFont{6}{7.2}{\rmdefault}{\mddefault}{\updefault}1}}}
\put(4876,-736){\makebox(0,0)[lb]{\smash{\SetFigFont{8}{9.6}{\rmdefault}{\mddefault}{\updefault}p}}}
\put(7726,-886){\makebox(0,0)[lb]{\smash{\SetFigFont{8}{9.6}{\rmdefault}{\mddefault}{\updefault}p}}}
\put(8326,-886){\makebox(0,0)[lb]{\smash{\SetFigFont{8}{9.6}{\rmdefault}{\mddefault}{\updefault}p}}}
\put(10801,-811){\makebox(0,0)[lb]{\smash{\SetFigFont{8}{9.6}{\rmdefault}{\mddefault}{\updefault}p}}}
\put(11101,239){\makebox(0,0)[lb]{\smash{\SetFigFont{8}{9.6}{\rmdefault}{\mddefault}{\updefault}1}}}
\put(11776,239){\makebox(0,0)[lb]{\smash{\SetFigFont{8}{9.6}{\rmdefault}{\mddefault}{\updefault}1}}}
\put(7951,239){\makebox(0,0)[lb]{\smash{\SetFigFont{8}{9.6}{\rmdefault}{\mddefault}{\updefault}1}}}
\put(8476,-1036){\makebox(0,0)[lb]{\smash{\SetFigFont{6}{7.2}{\rmdefault}{\mddefault}{\updefault}i}}}
\put(10951,-961){\makebox(0,0)[lb]{\smash{\SetFigFont{6}{7.2}{\rmdefault}{\mddefault}{\updefault}j}}}
\put(7801,-1036){\makebox(0,0)[lb]{\smash{\SetFigFont{6}{7.2}{\rmdefault}{\mddefault}{\updefault}i+1}}}
\put(7351,239){\makebox(0,0)[lb]{\smash{\SetFigFont{8}{9.6}{\rmdefault}{\mddefault}{\updefault}1}}}
\put(8701,314){\makebox(0,0)[lb]{\smash{\SetFigFont{8}{9.6}{\rmdefault}{\mddefault}{\updefault}1}}}
\put(9526,239){\makebox(0,0)[lb]{\smash{\SetFigFont{8}{9.6}{\rmdefault}{\mddefault}{\updefault}1}}}
\put(10351,314){\makebox(0,0)[lb]{\smash{\SetFigFont{8}{9.6}{\rmdefault}{\mddefault}{\updefault}1}}}
\end{picture}
\end{center}

Remark that the $T_{\gamma_i}$'s commute with all the other 
$T_{\gamma_{jk}}$ because their support curves can be made disjoint.
These formulas make up a representation of the pure braid 
group $P_n$, which extends to the whole braid group in the 
obvious manner: consider that the $i$-th and $i+1$-th legs are only
half-twisted. This defines the action of the $i$-th generator $g_i$ of
the braid group $B_n$. In fact, looking at the generators $A_{ij}$ of 
$P_n$ as elements of $B_n$, their action on $V$ consists in twisting 
the corresponding strands of $\Gamma'$, modulo Reidemester moves in
plane. On the other hand the fact that we obtained a representation of
$B_n$ is checked the same manner: the relation
$g_ig_{i+1}g_i=g_{i+1}g_ig_{i+1}$ translates into the third 
Reidemester move, which is obviously satisfied in the skein module. 
We continue to denote by $\rho$ the representation of $B_n$ on $V$.
In fact this should enter in the computation of the action  
of elements in the larger group ${\cal M}(F)$, so actually is ``part''
of $\rho$ but for bigger genus. 

The main ingredient in the proof of proposition \ref{pr} is:
\newtheorem{prop1}[prop]{Proposition}
\begin{prop1}\label{pr2}
The representations 
\[ (-A^{-1})\rho:B_3\longrightarrow End(V(3,1))\]
and respectively 
\[ (-A^{-1})\rho:B_4\longrightarrow End(V(4,2))\]
factor throughout the Temperley-Lieb algebras 
$A_{\beta,3}$ (and respectively $A_{\beta,4}$) 
where $\beta=2+q+q^{-1}$ and the parameter $q=A^{-4}$.
\end{prop1}
Let us explain what is meant by $c\rho$ where $c\in {\bf C}$ is a constant. 
This is another representation of the braid group, which is defined 
on the generators by $(c\rho)(g_i)=c \rho(g_i)$. Therefore if $w$ is a
word in the generators $g_i$ we have 
$(c\rho)(w)=c^{\mid w\mid} \rho(w)$, where $\mid w\mid$ is the sum of 
exponents appearing in the word $w$. Since the braid  relations are
homogeneous in the generators $g_i$ this is well-defined. 

We will prove that the image   of the braid group is infinite 
in $A_{\beta,n}$. 
This was done by Jones in \cite{Jones2}
for one value of $A$ (and a slightly different context), 
but the proof extends to an arbitrary
primitive root of unity. We can be more precise. Let us consider 
$M_2({\bf C})$ as the image of $A_{\beta,3}$ by the natural 
projection, and the factor  $M_3({\bf C})$ as the quotient of
$A_{\beta,4}$,
when both are semi-simple.   

\newtheorem{prop2}[prop]{Proposition}
\begin{prop2}\label{pr3}
The image of  $B_3$ in $M_2({\bf C})$ (via  $A_{\beta,3}$) 
is infinite provided that $r\neq 2,3,4,6,10$. 
For $r=10$ the image of  $B_4$ in $M_3({\bf C})$ (via  $A_{\beta,4}$) 
is infinite. 
\end{prop2}

Now the algebra $A_{\beta,n}$ ($n= 3,4$) is semi-simple 
for $r$ outside the excepted range $r\neq 3,4$ and  
$A_{\beta,3}= M_2({\bf C})\oplus {\bf C}$. 
The representation $\tilde{\rho}$ is not abelian and 
then the induced map 
 $A_{\beta,3}\longrightarrow M_2({\bf C})=End(V(3,1))$ 
must be  the canonical projection (up to an automorphism). 
Thus the image of $B_3$ by $\rho$ should be infinite. 
For $r=10$ we work within  $B_4$ and $A_{\beta,4}$ and we are 
forced then to restrict to  
$g\geq 3$. Again $\tilde{\rho}$ is not abelian and we will see 
in the next section that it is irreducible. Therefore the induced map 
$A_{\beta,4}\longrightarrow M_3({\bf C})=End(V(4,2))$ 
is again the projection on the corresponding factor.   
This establishes the proposition \ref{pr}, because $P_n$ is of finite
index in $B_n$. 

Eventually recall that $\rho$, at the mapping class group level, 
is only a projective representation, and it  can be also
considered  as a representation of a finite extension (depending on the level) 
of ${\cal M}_g$. The image group stays then in the unitary 
group modulo roots of unity of order $4r$ and thus  the claim of theorem 1
(concerning the BHMV-version) is proved. The present  proof (see the section 3.3) shows actually that
the image in the projective unitary group (i.e. modulo 
$U(1)$) is also infinite. $\Box$

\subsection{Proof of proposition \ref{pr2}}
The TQFT considered here is the one constructed in \cite{BHMV2}, for
$A$ a primitive $2r$-th root of unity. We  suppose for 
simplicity  that $r$ is
even, hence we are working with the $SU(2)$-TQFT. The same  arguments
hold verbatim  for the 
representation associated to the $SO(3)$-TQFT, 
with only minor modifications in the range of colors.  
\newtheorem{lem1}{Lemma}[subsection]
\begin{lem1}
A basis for $V(n,m)$ is provided by the labeled  graphs $L({\bf p})$ below,

\vspace{0.5cm}
\begin{center}
\setlength{\unitlength}{0.00025000in}%
\begingroup\makeatletter\ifx\SetFigFont\undefined%
\gdef\SetFigFont#1#2#3#4#5{%
  \reset@font\fontsize{#1}{#2pt}%
  \fontfamily{#3}\fontseries{#4}\fontshape{#5}%
  \selectfont}%
\fi\endgroup%
\begin{picture}(5274,1395)(1264,-1219)
\thicklines
\put(1276,-511){\line( 1, 0){5250}}
\put(6526,-511){\line( 0, 1){  0}}
\put(2026,-511){\line( 0, 1){675}}
\put(2851,-511){\line( 0, 1){675}}
\put(3676,-511){\line( 0, 1){675}}
\put(5701,-511){\line( 0, 1){675}}
\put(4951,-511){\line( 0, 1){675}}
\put(2326,-961){\makebox(0,0)[lb]{\smash{\SetFigFont{6}{7.2}{\rmdefault}{\mddefault}{\updefault}p}}}
\put(3151,-961){\makebox(0,0)[lb]{\smash{\SetFigFont{6}{7.2}{\rmdefault}{\mddefault}{\updefault}p}}}
\put(5251,-961){\makebox(0,0)[lb]{\smash{\SetFigFont{6}{7.2}{\rmdefault}{\mddefault}{\updefault}p}}}
\put(6151,-961){\makebox(0,0)[lb]{\smash{\SetFigFont{6}{7.2}{\rmdefault}{\mddefault}{\updefault}0}}}
\put(1426,-961){\makebox(0,0)[lb]{\smash{\SetFigFont{6}{7.2}{\rmdefault}{\mddefault}{\updefault}m}}}
\put(1801, 14){\makebox(0,0)[lb]{\smash{\SetFigFont{6}{7.2}{\rmdefault}{\mddefault}{\updefault}1}}}
\put(2626, 14){\makebox(0,0)[lb]{\smash{\SetFigFont{6}{7.2}{\rmdefault}{\mddefault}{\updefault}1}}}
\put(3451, 14){\makebox(0,0)[lb]{\smash{\SetFigFont{6}{7.2}{\rmdefault}{\mddefault}{\updefault}1}}}
\put(4726, 14){\makebox(0,0)[lb]{\smash{\SetFigFont{6}{7.2}{\rmdefault}{\mddefault}{\updefault}1}}}
\put(5476, 14){\makebox(0,0)[lb]{\smash{\SetFigFont{6}{7.2}{\rmdefault}{\mddefault}{\updefault}1}}}
\put(3301,-1186){\makebox(0,0)[lb]{\smash{\SetFigFont{5}{6.0}{\rmdefault}{\mddefault}{\updefault}n-2}}}
\put(2401,-1186){\makebox(0,0)[lb]{\smash{\SetFigFont{5}{6.0}{\rmdefault}{\mddefault}{\updefault}n-1}}}
\put(5326,-1186){\makebox(0,0)[lb]{\smash{\SetFigFont{5}{6.0}{\rmdefault}{\mddefault}{\updefault}1}}}
\end{picture}
\end{center}

\noindent  whose labels are in one-to-one correspondence with 
\[B(V)=\{{\bf p}=(p_0,p_1,....,p_{n}); p_i\in {\bf Z}_+, p_0=0, p_n=m,
p_i\leq 2r-2, \mid p_i-p_{i+1}\mid = 1, i=0,...,n\}.\]

\end{lem1}
{\em Proof:} It follows immediately from the admissibility conditions 
on the triples $(p_i,p_{i+1},1)$. $\Box$ \\
Then the computation of $\rho(g_i)$ is reduced to that of $g_i L({\bf
  p})$, in the skein module. Observe now that the only values of the
  labels $p_j$ which may change when $g_i$ is applied are $p_{i-1},
  p_i, p_{i+1}$. This will be also seen during the explicit
  computation.

Actually we have to  compute 
\setlength{\unitlength}{0.00050000in}%
\begingroup\makeatletter\ifx\SetFigFont\undefined%
\gdef\SetFigFont#1#2#3#4#5{%
  \reset@font\fontsize{#1}{#2pt}%
  \fontfamily{#3}\fontseries{#4}\fontshape{#5}%
  \selectfont}%
\fi\endgroup%
\begin{picture}(1074,864)(1564,-1063)
\thicklines
\put(1876,-736){\line( 3, 4){387}}
\multiput(2251,-736)(-7.89474,11.84211){20}{\makebox(11.1111,16.6667){\SetFigFont{7}{8.4}{\rmdefault}{\mddefault}{\updefault}.}}
\multiput(2026,-436)(-7.89474,11.84211){20}{\makebox(11.1111,16.6667){\SetFigFont{7}{8.4}{\rmdefault}{\mddefault}{\updefault}.}}
\put(1576,-736){\line( 1, 0){1050}}
\put(2026,-1036){\makebox(0,0)[lb]{\smash{\SetFigFont{7}{8.4}{\rmdefault}{\mddefault}{\updefault}b}}}
\put(2401,-1036){\makebox(0,0)[lb]{\smash{\SetFigFont{7}{8.4}{\rmdefault}{\mddefault}{\updefault}c}}}
\put(1651,-1036){\makebox(0,0)[lb]{\smash{\SetFigFont{7}{8.4}{\rmdefault}{\mddefault}{\updefault}a}}}
\put(1726,-436){\makebox(0,0)[lb]{\smash{\SetFigFont{7}{8.4}{\rmdefault}{\mddefault}{\updefault}1}}}
\put(2326,-436){\makebox(0,0)[lb]{\smash{\SetFigFont{7}{8.4}{\rmdefault}{\mddefault}{\updefault}1}}}
\end{picture}
and according to the lemma it suffices to consider that 
 $\mid a-b\mid=\mid b-c\mid =1$. 
\newtheorem{lem2}[lem1]{Lemma}
\begin{lem2}\label{lem2}
Suppose that $\mid a-b\mid=\mid b-c\mid =1$ and $\mid a-c\mid =2$. 
Then 

\vspace{0.5cm}
\hspace{3cm} 
\parbox{2cm}{\setlength{\unitlength}{0.00050000in}%
\begingroup\makeatletter\ifx\SetFigFont\undefined%
\gdef\SetFigFont#1#2#3#4#5{%
  \reset@font\fontsize{#1}{#2pt}%
  \fontfamily{#3}\fontseries{#4}\fontshape{#5}%
  \selectfont}%
\fi\endgroup%
\begin{picture}(1074,864)(1564,-1063)
\thicklines
\put(1876,-736){\line( 3, 4){387}}
\multiput(2251,-736)(-7.89474,11.84211){20}{\makebox(11.1111,16.6667){\SetFigFont{7}{8.4}{\rmdefault}{\mddefault}{\updefault}.}}
\multiput(2026,-436)(-7.89474,11.84211){20}{\makebox(11.1111,16.6667){\SetFigFont{7}{8.4}{\rmdefault}{\mddefault}{\updefault}.}}
\put(1576,-736){\line( 1, 0){1050}}
\put(2026,-1036){\makebox(0,0)[lb]{\smash{\SetFigFont{7}{8.4}{\rmdefault}{\mddefault}{\updefault}b}}}
\put(2401,-1036){\makebox(0,0)[lb]{\smash{\SetFigFont{7}{8.4}{\rmdefault}{\mddefault}{\updefault}c}}}
\put(1651,-1036){\makebox(0,0)[lb]{\smash{\SetFigFont{7}{8.4}{\rmdefault}{\mddefault}{\updefault}a}}}
\put(1726,-436){\makebox(0,0)[lb]{\smash{\SetFigFont{7}{8.4}{\rmdefault}{\mddefault}{\updefault}1}}}
\put(2326,-436){\makebox(0,0)[lb]{\smash{\SetFigFont{7}{8.4}{\rmdefault}{\mddefault}{\updefault}1}}}
\end{picture}}
 $\;\; = \; A\;\;$ 
\parbox{2cm}{\setlength{\unitlength}{0.00050000in}%
\begingroup\makeatletter\ifx\SetFigFont\undefined%
\gdef\SetFigFont#1#2#3#4#5{%
  \reset@font\fontsize{#1}{#2pt}%
  \fontfamily{#3}\fontseries{#4}\fontshape{#5}%
  \selectfont}%
\fi\endgroup%
\begin{picture}(1074,714)(1564,-1063)
\thicklines
\put(1576,-736){\line( 1, 0){1050}}
\put(1876,-736){\line( 0, 1){375}}
\put(2251,-736){\line( 0, 1){375}}
\put(2026,-1036){\makebox(0,0)[lb]{\smash{\SetFigFont{7}{8.4}{\rmdefault}{\mddefault}{\updefault}b}}}
\put(2401,-1036){\makebox(0,0)[lb]{\smash{\SetFigFont{7}{8.4}{\rmdefault}{\mddefault}{\updefault}c}}}
\put(1651,-1036){\makebox(0,0)[lb]{\smash{\SetFigFont{7}{8.4}{\rmdefault}{\mddefault}{\updefault}a}}}
\end{picture}}

\end{lem2}
{\em Proof:} Suppose for simplicity that $c=a+2, b=a+1$. Then according to
\cite{MV} and the  section 2.3 we have 

\vspace{0.5cm}
\hspace{3cm} 
\parbox{2cm}{\setlength{\unitlength}{0.00050000in}%
\begingroup\makeatletter\ifx\SetFigFont\undefined%
\gdef\SetFigFont#1#2#3#4#5{%
  \reset@font\fontsize{#1}{#2pt}%
  \fontfamily{#3}\fontseries{#4}\fontshape{#5}%
  \selectfont}%
\fi\endgroup%
\begin{picture}(1074,864)(1564,-1063)
\thicklines
\multiput(2251,-736)(-7.89474,11.84211){20}{\makebox(11.1111,16.6667){\SetFigFont{7}{8.4}{\rmdefault}{\mddefault}{\updefault}.}}
\put(1576,-736){\line( 1, 0){1050}}
\put(1876,-736){\line( 3, 4){387}}
\multiput(1876,-211)(7.89474,-11.84211){20}{\makebox(11.1111,16.6667){\SetFigFont{7}{8.4}{\rmdefault}{\mddefault}{\updefault}.}}
\put(1876,-1036){\makebox(0,0)[lb]{\smash{\SetFigFont{7}{8.4}{\rmdefault}{\mddefault}{\updefault}a+1}}}
\put(2401,-1036){\makebox(0,0)[lb]{\smash{\SetFigFont{7}{8.4}{\rmdefault}{\mddefault}{\updefault}a+2}}}
\put(1576,-1036){\makebox(0,0)[lb]{\smash{\SetFigFont{7}{8.4}{\rmdefault}{\mddefault}{\updefault}a}}}
\put(1801,-436){\makebox(0,0)[lb]{\smash{\SetFigFont{7}{8.4}{\rmdefault}{\mddefault}{\updefault}1}}}
\put(2326,-436){\makebox(0,0)[lb]{\smash{\SetFigFont{7}{8.4}{\rmdefault}{\mddefault}{\updefault}1}}}
\end{picture}} 
\parbox{5cm}{$\; = \; \; \sum_k \delta(k;1,1)\frac{<k>}{<1,1,k>}\;   $}
\parbox{2cm}{\setlength{\unitlength}{0.00050000in}%
\begingroup\makeatletter\ifx\SetFigFont\undefined%
\gdef\SetFigFont#1#2#3#4#5{%
  \reset@font\fontsize{#1}{#2pt}%
  \fontfamily{#3}\fontseries{#4}\fontshape{#5}%
  \selectfont}%
\fi\endgroup%
\begin{picture}(1074,1110)(1564,-1063)
\thicklines
\put(1576,-736){\line( 1, 0){1050}}
\multiput(1876,-736)(10.22727,10.22727){23}{\makebox(11.1111,16.6667){\SetFigFont{7}{8.4}{\rmdefault}{\mddefault}{\updefault}.}}
\multiput(2101,-511)(10.22727,-10.22727){23}{\makebox(11.1111,16.6667){\SetFigFont{7}{8.4}{\rmdefault}{\mddefault}{\updefault}.}}
\put(2101,-511){\line( 0, 1){300}}
\multiput(2101,-211)(10.22727,10.22727){23}{\makebox(11.1111,16.6667){\SetFigFont{7}{8.4}{\rmdefault}{\mddefault}{\updefault}.}}
\multiput(1876, 14)(10.22727,-10.22727){23}{\makebox(11.1111,16.6667){\SetFigFont{7}{8.4}{\rmdefault}{\mddefault}{\updefault}.}}
\put(1876,-1036){\makebox(0,0)[lb]{\smash{\SetFigFont{7}{8.4}{\rmdefault}{\mddefault}{\updefault}a+1}}}
\put(2401,-1036){\makebox(0,0)[lb]{\smash{\SetFigFont{7}{8.4}{\rmdefault}{\mddefault}{\updefault}a+2}}}
\put(1576,-1036){\makebox(0,0)[lb]{\smash{\SetFigFont{7}{8.4}{\rmdefault}{\mddefault}{\updefault}a}}}
\put(2176,-436){\makebox(0,0)[lb]{\smash{\SetFigFont{7}{8.4}{\rmdefault}{\mddefault}{\updefault}k}}}
\put(2401,-61){\makebox(0,0)[lb]{\smash{\SetFigFont{7}{8.4}{\rmdefault}{\mddefault}{\updefault}1}}}
\put(1726,-61){\makebox(0,0)[lb]{\smash{\SetFigFont{7}{8.4}{\rmdefault}{\mddefault}{\updefault}1}}}
\put(1876,-586){\makebox(0,0)[lb]{\smash{\SetFigFont{7}{8.4}{\rmdefault}{\mddefault}{\updefault}1}}}
\put(2326,-586){\makebox(0,0)[lb]{\smash{\SetFigFont{7}{8.4}{\rmdefault}{\mddefault}{\updefault}1}}}
\end{picture}}
\vspace{0.5cm}

\noindent Therefore the triple $(1,1,k)$ has to be  admissible, so that 
$k\in\{0,2\}$. Also the open tetrahedron 
\begin{center}
\setlength{\unitlength}{0.00050000in}%
\begingroup\makeatletter\ifx\SetFigFont\undefined%
\gdef\SetFigFont#1#2#3#4#5{%
  \reset@font\fontsize{#1}{#2pt}%
  \fontfamily{#3}\fontseries{#4}\fontshape{#5}%
  \selectfont}%
\fi\endgroup%
\begin{picture}(1074,885)(1564,-1063)
\thicklines
\put(1576,-736){\line( 1, 0){1050}}
\multiput(1876,-736)(10.22727,10.22727){23}{\makebox(11.1111,16.6667){\SetFigFont{7}{8.4}{\rmdefault}{\mddefault}{\updefault}.}}
\multiput(2101,-511)(10.22727,-10.22727){23}{\makebox(11.1111,16.6667){\SetFigFont{7}{8.4}{\rmdefault}{\mddefault}{\updefault}.}}
\put(2101,-511){\line( 0, 1){300}}
\put(1876,-1036){\makebox(0,0)[lb]{\smash{\SetFigFont{7}{8.4}{\rmdefault}{\mddefault}{\updefault}a+1}}}
\put(2401,-1036){\makebox(0,0)[lb]{\smash{\SetFigFont{7}{8.4}{\rmdefault}{\mddefault}{\updefault}a+2}}}
\put(1576,-1036){\makebox(0,0)[lb]{\smash{\SetFigFont{7}{8.4}{\rmdefault}{\mddefault}{\updefault}a}}}
\put(1876,-586){\makebox(0,0)[lb]{\smash{\SetFigFont{7}{8.4}{\rmdefault}{\mddefault}{\updefault}1}}}
\put(2326,-586){\makebox(0,0)[lb]{\smash{\SetFigFont{7}{8.4}{\rmdefault}{\mddefault}{\updefault}1}}}
\put(2176,-286){\makebox(0,0)[lb]{\smash{\SetFigFont{7}{8.4}{\rmdefault}{\mddefault}{\updefault}k}}}
\end{picture}  
\end{center}

\noindent vanishes if $(k,a,a+2)$ is not admissible, so that there is only one
possibility left, namely $k=2$. We get rid of the triangular face by using
the formula:

\vspace{0.5cm}
\hspace{3cm} 
\parbox{2cm}{\setlength{\unitlength}{0.00050000in}%
\begingroup\makeatletter\ifx\SetFigFont\undefined%
\gdef\SetFigFont#1#2#3#4#5{%
  \reset@font\fontsize{#1}{#2pt}%
  \fontfamily{#3}\fontseries{#4}\fontshape{#5}%
  \selectfont}%
\fi\endgroup%
\begin{picture}(1074,1110)(1564,-1063)
\thicklines
\put(1576,-736){\line( 1, 0){1050}}
\multiput(1876,-736)(10.22727,10.22727){23}{\makebox(11.1111,16.6667){\SetFigFont{7}{8.4}{\rmdefault}{\mddefault}{\updefault}.}}
\multiput(2101,-511)(10.22727,-10.22727){23}{\makebox(11.1111,16.6667){\SetFigFont{7}{8.4}{\rmdefault}{\mddefault}{\updefault}.}}
\put(2101,-511){\line( 0, 1){300}}
\multiput(2101,-211)(10.22727,10.22727){23}{\makebox(11.1111,16.6667){\SetFigFont{7}{8.4}{\rmdefault}{\mddefault}{\updefault}.}}
\multiput(1876, 14)(10.22727,-10.22727){23}{\makebox(11.1111,16.6667){\SetFigFont{7}{8.4}{\rmdefault}{\mddefault}{\updefault}.}}
\put(1876,-1036){\makebox(0,0)[lb]{\smash{\SetFigFont{7}{8.4}{\rmdefault}{\mddefault}{\updefault}a+1}}}
\put(2401,-1036){\makebox(0,0)[lb]{\smash{\SetFigFont{7}{8.4}{\rmdefault}{\mddefault}{\updefault}a+2}}}
\put(1576,-1036){\makebox(0,0)[lb]{\smash{\SetFigFont{7}{8.4}{\rmdefault}{\mddefault}{\updefault}a}}}
\put(2401,-61){\makebox(0,0)[lb]{\smash{\SetFigFont{7}{8.4}{\rmdefault}{\mddefault}{\updefault}1}}}
\put(1726,-61){\makebox(0,0)[lb]{\smash{\SetFigFont{7}{8.4}{\rmdefault}{\mddefault}{\updefault}1}}}
\put(1876,-586){\makebox(0,0)[lb]{\smash{\SetFigFont{7}{8.4}{\rmdefault}{\mddefault}{\updefault}1}}}
\put(2326,-586){\makebox(0,0)[lb]{\smash{\SetFigFont{7}{8.4}{\rmdefault}{\mddefault}{\updefault}1}}}
\put(2176,-361){\makebox(0,0)[lb]{\smash{\SetFigFont{7}{8.4}{\rmdefault}{\mddefault}{\updefault}2}}}
\end{picture}} 
\parbox{5cm}{$\; = \; \; \frac{\left < 
                \begin{array}{ccc}
                 2 & 1 & 1 \\
                a+1 & a+2 & a
                 \end{array} \right >}{<2,a,a+2>}\;   $}
\parbox{2cm}{\setlength{\unitlength}{0.00050000in}%
\begingroup\makeatletter\ifx\SetFigFont\undefined%
\gdef\SetFigFont#1#2#3#4#5{%
  \reset@font\fontsize{#1}{#2pt}%
  \fontfamily{#3}\fontseries{#4}\fontshape{#5}%
  \selectfont}%
\fi\endgroup%
\begin{picture}(762,864)(1726,-1063)
\thicklines
\put(2101,-811){\makebox(11.1111,16.6667){\SetFigFont{10}{12}{\rmdefault}{\mddefault}{\updefault}.}}
\put(2101,-811){\line( 0, 1){450}}
\multiput(2101,-361)(-11.84211,7.89474){20}{\makebox(11.1111,16.6667){\SetFigFont{7}{8.4}{\rmdefault}{\mddefault}{\updefault}.}}
\put(1801,-811){\line( 1, 0){675}}
\multiput(2101,-361)(11.84211,7.89474){20}{\makebox(11.1111,16.6667){\SetFigFont{7}{8.4}{\rmdefault}{\mddefault}{\updefault}.}}
\put(1876,-1036){\makebox(0,0)[lb]{\smash{\SetFigFont{7}{8.4}{\rmdefault}{\mddefault}{\updefault}a}}}
\put(2176,-1036){\makebox(0,0)[lb]{\smash{\SetFigFont{7}{8.4}{\rmdefault}{\mddefault}{\updefault}a+2}}}
\put(2176,-661){\makebox(0,0)[lb]{\smash{\SetFigFont{7}{8.4}{\rmdefault}{\mddefault}{\updefault}2}}}
\put(1726,-361){\makebox(0,0)[lb]{\smash{\SetFigFont{7}{8.4}{\rmdefault}{\mddefault}{\updefault}1}}}
\put(2401,-361){\makebox(0,0)[lb]{\smash{\SetFigFont{7}{8.4}{\rmdefault}{\mddefault}{\updefault}1}}}
\end{picture}}
\vspace{0.5cm}

\noindent
We eventually  perform a fusion in order to express the 
right hand side in terms of the the usual basis
 $L({\bf p})$. The formula of the fusing is:

\vspace{0.5cm}
\hspace{3cm} 
\parbox{2cm}{\setlength{\unitlength}{0.00050000in}%
\begingroup\makeatletter\ifx\SetFigFont\undefined%
\gdef\SetFigFont#1#2#3#4#5{%
  \reset@font\fontsize{#1}{#2pt}%
  \fontfamily{#3}\fontseries{#4}\fontshape{#5}%
  \selectfont}%
\fi\endgroup%
\begin{picture}(1074,714)(1564,-1063)
\thicklines
\put(1576,-736){\line( 1, 0){1050}}
\put(1876,-736){\line( 0, 1){375}}
\put(2251,-736){\line( 0, 1){375}}
\put(1651,-1036){\makebox(0,0)[lb]{\smash{\SetFigFont{7}{8.4}{\rmdefault}{\mddefault}{\updefault}a}}}
\put(2326,-1036){\makebox(0,0)[lb]{\smash{\SetFigFont{7}{8.4}{\rmdefault}{\mddefault}{\updefault}a+2}}}
\put(1951,-1036){\makebox(0,0)[lb]{\smash{\SetFigFont{7}{8.4}{\rmdefault}{\mddefault}{\updefault}a+1}}}
\put(1726,-511){\makebox(0,0)[lb]{\smash{\SetFigFont{7}{8.4}{\rmdefault}{\mddefault}{\updefault}1}}}
\put(2326,-511){\makebox(0,0)[lb]{\smash{\SetFigFont{7}{8.4}{\rmdefault}{\mddefault}{\updefault}1}}}
\end{picture}} 
\parbox{5cm}{$\; = \; \; \left \{ 
                \begin{array}{ccc}
                 a & 1 & 2 \\
                 1 & a+2 & a+1
                 \end{array} \right \}\;   $}
\parbox{2cm}{\setlength{\unitlength}{0.00050000in}%
\begingroup\makeatletter\ifx\SetFigFont\undefined%
\gdef\SetFigFont#1#2#3#4#5{%
  \reset@font\fontsize{#1}{#2pt}%
  \fontfamily{#3}\fontseries{#4}\fontshape{#5}%
  \selectfont}%
\fi\endgroup%
\begin{picture}(762,864)(1726,-1063)
\thicklines
\put(2101,-811){\makebox(11.1111,16.6667){\SetFigFont{10}{12}{\rmdefault}{\mddefault}{\updefault}.}}
\put(2101,-811){\line( 0, 1){450}}
\multiput(2101,-361)(-11.84211,7.89474){20}{\makebox(11.1111,16.6667){\SetFigFont{7}{8.4}{\rmdefault}{\mddefault}{\updefault}.}}
\put(1801,-811){\line( 1, 0){675}}
\multiput(2101,-361)(11.84211,7.89474){20}{\makebox(11.1111,16.6667){\SetFigFont{7}{8.4}{\rmdefault}{\mddefault}{\updefault}.}}
\put(1876,-1036){\makebox(0,0)[lb]{\smash{\SetFigFont{7}{8.4}{\rmdefault}{\mddefault}{\updefault}a}}}
\put(2176,-1036){\makebox(0,0)[lb]{\smash{\SetFigFont{7}{8.4}{\rmdefault}{\mddefault}{\updefault}a+2}}}
\put(2176,-661){\makebox(0,0)[lb]{\smash{\SetFigFont{7}{8.4}{\rmdefault}{\mddefault}{\updefault}2}}}
\put(1726,-361){\makebox(0,0)[lb]{\smash{\SetFigFont{7}{8.4}{\rmdefault}{\mddefault}{\updefault}1}}}
\put(2401,-361){\makebox(0,0)[lb]{\smash{\SetFigFont{7}{8.4}{\rmdefault}{\mddefault}{\updefault}1}}}
\end{picture}}
\vspace{0.5cm} \\
where the quantum 6-j symbol involved, namely 
$\left\{\begin{array}{ccc}
a & 1 & 2 \\
1 & a+2 & a+1 
\end{array}\right\}$, it turns out to be equal to $1$, for all $a$. 
This implies that

\vspace{0.5cm}
\hspace{3cm} 
\parbox{2cm}{\setlength{\unitlength}{0.00050000in}%
\begingroup\makeatletter\ifx\SetFigFont\undefined%
\gdef\SetFigFont#1#2#3#4#5{%
  \reset@font\fontsize{#1}{#2pt}%
  \fontfamily{#3}\fontseries{#4}\fontshape{#5}%
  \selectfont}%
\fi\endgroup%
\begin{picture}(1074,1110)(1564,-1063)
\thicklines
\put(1576,-736){\line( 1, 0){1050}}
\multiput(1876,-736)(10.22727,10.22727){23}{\makebox(11.1111,16.6667){\SetFigFont{7}{8.4}{\rmdefault}{\mddefault}{\updefault}.}}
\multiput(2101,-511)(10.22727,-10.22727){23}{\makebox(11.1111,16.6667){\SetFigFont{7}{8.4}{\rmdefault}{\mddefault}{\updefault}.}}
\put(2101,-511){\line( 0, 1){300}}
\multiput(2101,-211)(10.22727,10.22727){23}{\makebox(11.1111,16.6667){\SetFigFont{7}{8.4}{\rmdefault}{\mddefault}{\updefault}.}}
\multiput(1876, 14)(10.22727,-10.22727){23}{\makebox(11.1111,16.6667){\SetFigFont{7}{8.4}{\rmdefault}{\mddefault}{\updefault}.}}
\put(1876,-1036){\makebox(0,0)[lb]{\smash{\SetFigFont{7}{8.4}{\rmdefault}{\mddefault}{\updefault}a+1}}}
\put(2401,-1036){\makebox(0,0)[lb]{\smash{\SetFigFont{7}{8.4}{\rmdefault}{\mddefault}{\updefault}a+2}}}
\put(1576,-1036){\makebox(0,0)[lb]{\smash{\SetFigFont{7}{8.4}{\rmdefault}{\mddefault}{\updefault}a}}}
\put(2401,-61){\makebox(0,0)[lb]{\smash{\SetFigFont{7}{8.4}{\rmdefault}{\mddefault}{\updefault}1}}}
\put(1726,-61){\makebox(0,0)[lb]{\smash{\SetFigFont{7}{8.4}{\rmdefault}{\mddefault}{\updefault}1}}}
\put(1876,-586){\makebox(0,0)[lb]{\smash{\SetFigFont{7}{8.4}{\rmdefault}{\mddefault}{\updefault}1}}}
\put(2326,-586){\makebox(0,0)[lb]{\smash{\SetFigFont{7}{8.4}{\rmdefault}{\mddefault}{\updefault}1}}}
\put(2176,-361){\makebox(0,0)[lb]{\smash{\SetFigFont{7}{8.4}{\rmdefault}{\mddefault}{\updefault}2}}}
\end{picture}} 
\parbox{2cm}{$\; = \; $}
\parbox{2cm}{\setlength{\unitlength}{0.00050000in}%
\begingroup\makeatletter\ifx\SetFigFont\undefined%
\gdef\SetFigFont#1#2#3#4#5{%
  \reset@font\fontsize{#1}{#2pt}%
  \fontfamily{#3}\fontseries{#4}\fontshape{#5}%
  \selectfont}%
\fi\endgroup%
\begin{picture}(1074,714)(1564,-1063)
\thicklines
\put(1576,-736){\line( 1, 0){1050}}
\put(1876,-736){\line( 0, 1){375}}
\put(2251,-736){\line( 0, 1){375}}
\put(1651,-1036){\makebox(0,0)[lb]{\smash{\SetFigFont{7}{8.4}{\rmdefault}{\mddefault}{\updefault}a}}}
\put(2326,-1036){\makebox(0,0)[lb]{\smash{\SetFigFont{7}{8.4}{\rmdefault}{\mddefault}{\updefault}a+2}}}
\put(1951,-1036){\makebox(0,0)[lb]{\smash{\SetFigFont{7}{8.4}{\rmdefault}{\mddefault}{\updefault}a+1}}}
\put(1726,-511){\makebox(0,0)[lb]{\smash{\SetFigFont{7}{8.4}{\rmdefault}{\mddefault}{\updefault}1}}}
\put(2326,-511){\makebox(0,0)[lb]{\smash{\SetFigFont{7}{8.4}{\rmdefault}{\mddefault}{\updefault}1}}}
\end{picture}}
\vspace{0.5cm} \\
and the lemma follows. $\Box$ \\

\newtheorem{lem3}[lem1]{Lemma}
\begin{lem3}\label{lem3}
The following identities hold: 

\vspace{0.5cm}
\begin{center} 
\parbox{2cm}{\setlength{\unitlength}{0.00050000in}%
\begingroup\makeatletter\ifx\SetFigFont\undefined%
\gdef\SetFigFont#1#2#3#4#5{%
  \reset@font\fontsize{#1}{#2pt}%
  \fontfamily{#3}\fontseries{#4}\fontshape{#5}%
  \selectfont}%
\fi\endgroup%
\begin{picture}(1074,717)(1489,-991)
\thicklines
\put(1801,-736){\line( 1, 1){450}}
\multiput(2251,-736)(-10.00000,10.00000){16}{\makebox(11.1111,16.6667){\SetFigFont{7}{8.4}{\rmdefault}{\mddefault}{\updefault}.}}
\multiput(1951,-436)(-10.00000,10.00000){16}{\makebox(11.1111,16.6667){\SetFigFont{7}{8.4}{\rmdefault}{\mddefault}{\updefault}.}}
\put(1501,-736){\line( 1, 0){1050}}
\put(2476,-961){\makebox(0,0)[lb]{\smash{\SetFigFont{7}{8.4}{\rmdefault}{\mddefault}{\updefault}a}}}
\put(1876,-961){\makebox(0,0)[lb]{\smash{\SetFigFont{8}{9.6}{\rmdefault}{\mddefault}{\updefault}a-1}}}
\put(1576,-961){\makebox(0,0)[lb]{\smash{\SetFigFont{8}{9.6}{\rmdefault}{\mddefault}{\updefault}a}}}
\put(1651,-436){\makebox(0,0)[lb]{\smash{\SetFigFont{8}{9.6}{\rmdefault}{\mddefault}{\updefault}1}}}
\put(2326,-436){\makebox(0,0)[lb]{\smash{\SetFigFont{8}{9.6}{\rmdefault}{\mddefault}{\updefault}1}}}
\end{picture}} 
\parbox{3cm}{$\; = \; -A\frac{A^{-4}[a]-[a+2]}{[2][a+1]}$}
\parbox{2cm}{\setlength{\unitlength}{0.00050000in}%
\begingroup\makeatletter\ifx\SetFigFont\undefined%
\gdef\SetFigFont#1#2#3#4#5{%
  \reset@font\fontsize{#1}{#2pt}%
  \fontfamily{#3}\fontseries{#4}\fontshape{#5}%
  \selectfont}%
\fi\endgroup%
\begin{picture}(1074,714)(1564,-1063)
\thicklines
\put(1576,-736){\line( 1, 0){1050}}
\put(1876,-736){\line( 0, 1){375}}
\put(2251,-736){\line( 0, 1){375}}
\put(1651,-1036){\makebox(0,0)[lb]{\smash{\SetFigFont{7}{8.4}{\rmdefault}{\mddefault}{\updefault}a}}}
\put(1726,-511){\makebox(0,0)[lb]{\smash{\SetFigFont{7}{8.4}{\rmdefault}{\mddefault}{\updefault}1}}}
\put(2326,-511){\makebox(0,0)[lb]{\smash{\SetFigFont{7}{8.4}{\rmdefault}{\mddefault}{\updefault}1}}}
\put(1951,-1036){\makebox(0,0)[lb]{\smash{\SetFigFont{7}{8.4}{\rmdefault}{\mddefault}{\updefault}a-1}}}
\put(2476,-1036){\makebox(0,0)[lb]{\smash{\SetFigFont{7}{8.4}{\rmdefault}{\mddefault}{\updefault}a}}}
\end{picture}}
\parbox{3cm}{$\; + \; A\frac{(A^{-4}+1)[a+1]}{[a+2]+[a]}$}
\parbox{2cm}{\setlength{\unitlength}{0.00050000in}%
\begingroup\makeatletter\ifx\SetFigFont\undefined%
\gdef\SetFigFont#1#2#3#4#5{%
  \reset@font\fontsize{#1}{#2pt}%
  \fontfamily{#3}\fontseries{#4}\fontshape{#5}%
  \selectfont}%
\fi\endgroup%
\begin{picture}(1074,714)(1564,-1063)
\thicklines
\put(1576,-736){\line( 1, 0){1050}}
\put(1876,-736){\line( 0, 1){375}}
\put(2251,-736){\line( 0, 1){375}}
\put(1651,-1036){\makebox(0,0)[lb]{\smash{\SetFigFont{7}{8.4}{\rmdefault}{\mddefault}{\updefault}a}}}
\put(1951,-1036){\makebox(0,0)[lb]{\smash{\SetFigFont{7}{8.4}{\rmdefault}{\mddefault}{\updefault}a+1}}}
\put(1726,-511){\makebox(0,0)[lb]{\smash{\SetFigFont{7}{8.4}{\rmdefault}{\mddefault}{\updefault}1}}}
\put(2326,-511){\makebox(0,0)[lb]{\smash{\SetFigFont{7}{8.4}{\rmdefault}{\mddefault}{\updefault}1}}}
\put(2401,-1036){\makebox(0,0)[lb]{\smash{\SetFigFont{7}{8.4}{\rmdefault}{\mddefault}{\updefault}a}}}
\end{picture}} 
\end{center}
\vspace{0.5cm} 

\vspace{0.5cm}
\begin{center} 
\parbox{2cm}{\setlength{\unitlength}{0.00050000in}%
\begingroup\makeatletter\ifx\SetFigFont\undefined%
\gdef\SetFigFont#1#2#3#4#5{%
  \reset@font\fontsize{#1}{#2pt}%
  \fontfamily{#3}\fontseries{#4}\fontshape{#5}%
  \selectfont}%
\fi\endgroup%
\begin{picture}(1074,717)(1489,-991)
\thicklines
\put(1801,-736){\line( 1, 1){450}}
\multiput(2251,-736)(-10.00000,10.00000){16}{\makebox(11.1111,16.6667){\SetFigFont{7}{8.4}{\rmdefault}{\mddefault}{\updefault}.}}
\multiput(1951,-436)(-10.00000,10.00000){16}{\makebox(11.1111,16.6667){\SetFigFont{7}{8.4}{\rmdefault}{\mddefault}{\updefault}.}}
\put(1501,-736){\line( 1, 0){1050}}
\put(2476,-961){\makebox(0,0)[lb]{\smash{\SetFigFont{7}{8.4}{\rmdefault}{\mddefault}{\updefault}a}}}
\put(1651,-436){\makebox(0,0)[lb]{\smash{\SetFigFont{8}{9.6}{\rmdefault}{\mddefault}{\updefault}1}}}
\put(2326,-436){\makebox(0,0)[lb]{\smash{\SetFigFont{8}{9.6}{\rmdefault}{\mddefault}{\updefault}1}}}
\put(1501,-961){\makebox(0,0)[lb]{\smash{\SetFigFont{8}{9.6}{\rmdefault}{\mddefault}{\updefault}a}}}
\put(1876,-961){\makebox(0,0)[lb]{\smash{\SetFigFont{8}{9.6}{\rmdefault}{\mddefault}{\updefault}a+1}}}
\end{picture}} 
\parbox{3cm}{$\; = \; A\frac{[a][a+2]}{A^2[a+1]^2}$}
\parbox{2cm}{\setlength{\unitlength}{0.00050000in}%
\begingroup\makeatletter\ifx\SetFigFont\undefined%
\gdef\SetFigFont#1#2#3#4#5{%
  \reset@font\fontsize{#1}{#2pt}%
  \fontfamily{#3}\fontseries{#4}\fontshape{#5}%
  \selectfont}%
\fi\endgroup%
\begin{picture}(1074,714)(1564,-1063)
\thicklines
\put(1576,-736){\line( 1, 0){1050}}
\put(1876,-736){\line( 0, 1){375}}
\put(2251,-736){\line( 0, 1){375}}
\put(1651,-1036){\makebox(0,0)[lb]{\smash{\SetFigFont{7}{8.4}{\rmdefault}{\mddefault}{\updefault}a}}}
\put(1726,-511){\makebox(0,0)[lb]{\smash{\SetFigFont{7}{8.4}{\rmdefault}{\mddefault}{\updefault}1}}}
\put(2326,-511){\makebox(0,0)[lb]{\smash{\SetFigFont{7}{8.4}{\rmdefault}{\mddefault}{\updefault}1}}}
\put(1951,-1036){\makebox(0,0)[lb]{\smash{\SetFigFont{7}{8.4}{\rmdefault}{\mddefault}{\updefault}a-1}}}
\put(2476,-1036){\makebox(0,0)[lb]{\smash{\SetFigFont{7}{8.4}{\rmdefault}{\mddefault}{\updefault}a}}}
\end{picture}}
\parbox{3cm}{$\; - \; A\frac{A^{-4}[a+2]-[a]}{[2][a+1]}$}
\parbox{2cm}{\setlength{\unitlength}{0.00050000in}%
\begingroup\makeatletter\ifx\SetFigFont\undefined%
\gdef\SetFigFont#1#2#3#4#5{%
  \reset@font\fontsize{#1}{#2pt}%
  \fontfamily{#3}\fontseries{#4}\fontshape{#5}%
  \selectfont}%
\fi\endgroup%
\begin{picture}(1074,714)(1564,-1063)
\thicklines
\put(1576,-736){\line( 1, 0){1050}}
\put(1876,-736){\line( 0, 1){375}}
\put(2251,-736){\line( 0, 1){375}}
\put(1651,-1036){\makebox(0,0)[lb]{\smash{\SetFigFont{7}{8.4}{\rmdefault}{\mddefault}{\updefault}a}}}
\put(1951,-1036){\makebox(0,0)[lb]{\smash{\SetFigFont{7}{8.4}{\rmdefault}{\mddefault}{\updefault}a+1}}}
\put(1726,-511){\makebox(0,0)[lb]{\smash{\SetFigFont{7}{8.4}{\rmdefault}{\mddefault}{\updefault}1}}}
\put(2326,-511){\makebox(0,0)[lb]{\smash{\SetFigFont{7}{8.4}{\rmdefault}{\mddefault}{\updefault}1}}}
\put(2401,-1036){\makebox(0,0)[lb]{\smash{\SetFigFont{7}{8.4}{\rmdefault}{\mddefault}{\updefault}a}}}
\end{picture}} 
\end{center}
\vspace{0.5cm} 

\noindent Here is to be understood that for $a=2r-2$ we have 
\parbox{2cm}{\setlength{\unitlength}{0.00050000in}%
\begingroup\makeatletter\ifx\SetFigFont\undefined%
\gdef\SetFigFont#1#2#3#4#5{%
  \reset@font\fontsize{#1}{#2pt}%
  \fontfamily{#3}\fontseries{#4}\fontshape{#5}%
  \selectfont}%
\fi\endgroup%
\begin{picture}(1074,714)(1564,-1063)
\thicklines
\put(1576,-736){\line( 1, 0){1050}}
\put(1876,-736){\line( 0, 1){375}}
\put(2251,-736){\line( 0, 1){375}}
\put(1651,-1036){\makebox(0,0)[lb]{\smash{\SetFigFont{7}{8.4}{\rmdefault}{\mddefault}{\updefault}a}}}
\put(1951,-1036){\makebox(0,0)[lb]{\smash{\SetFigFont{7}{8.4}{\rmdefault}{\mddefault}{\updefault}a+1}}}
\put(1726,-511){\makebox(0,0)[lb]{\smash{\SetFigFont{7}{8.4}{\rmdefault}{\mddefault}{\updefault}1}}}
\put(2326,-511){\makebox(0,0)[lb]{\smash{\SetFigFont{7}{8.4}{\rmdefault}{\mddefault}{\updefault}1}}}
\put(2401,-1036){\makebox(0,0)[lb]{\smash{\SetFigFont{7}{8.4}{\rmdefault}{\mddefault}{\updefault}a}}}
\end{picture}} 
\parbox{1cm}{$\; = \; 0 \;$,}
and for $a=0$ the  equality 
\parbox{2cm}{\setlength{\unitlength}{0.00050000in}%
\begingroup\makeatletter\ifx\SetFigFont\undefined%
\gdef\SetFigFont#1#2#3#4#5{%
  \reset@font\fontsize{#1}{#2pt}%
  \fontfamily{#3}\fontseries{#4}\fontshape{#5}%
  \selectfont}%
\fi\endgroup%
\begin{picture}(1074,714)(1564,-1063)
\thicklines
\put(1576,-736){\line( 1, 0){1050}}
\put(1876,-736){\line( 0, 1){375}}
\put(2251,-736){\line( 0, 1){375}}
\put(1651,-1036){\makebox(0,0)[lb]{\smash{\SetFigFont{7}{8.4}{\rmdefault}{\mddefault}{\updefault}a}}}
\put(1726,-511){\makebox(0,0)[lb]{\smash{\SetFigFont{7}{8.4}{\rmdefault}{\mddefault}{\updefault}1}}}
\put(2326,-511){\makebox(0,0)[lb]{\smash{\SetFigFont{7}{8.4}{\rmdefault}{\mddefault}{\updefault}1}}}
\put(1951,-1036){\makebox(0,0)[lb]{\smash{\SetFigFont{7}{8.4}{\rmdefault}{\mddefault}{\updefault}a-1}}}
\put(2476,-1036){\makebox(0,0)[lb]{\smash{\SetFigFont{7}{8.4}{\rmdefault}{\mddefault}{\updefault}a}}}
\end{picture}} 
\parbox{1cm}{$\; = \; 0 \;$,}
holds.
\end{lem3}
{\em Proof:} We have, like in the previous lemma, the following formula:

\vspace{0.5cm}
\begin{center} 
\parbox{2cm}{\setlength{\unitlength}{0.00050000in}%
\begingroup\makeatletter\ifx\SetFigFont\undefined%
\gdef\SetFigFont#1#2#3#4#5{%
  \reset@font\fontsize{#1}{#2pt}%
  \fontfamily{#3}\fontseries{#4}\fontshape{#5}%
  \selectfont}%
\fi\endgroup%
\begin{picture}(1074,717)(1489,-991)
\thicklines
\put(1801,-736){\line( 1, 1){450}}
\multiput(2251,-736)(-10.00000,10.00000){16}{\makebox(11.1111,16.6667){\SetFigFont{7}{8.4}{\rmdefault}{\mddefault}{\updefault}.}}
\multiput(1951,-436)(-10.00000,10.00000){16}{\makebox(11.1111,16.6667){\SetFigFont{7}{8.4}{\rmdefault}{\mddefault}{\updefault}.}}
\put(1501,-736){\line( 1, 0){1050}}
\put(2476,-961){\makebox(0,0)[lb]{\smash{\SetFigFont{7}{8.4}{\rmdefault}{\mddefault}{\updefault}a}}}
\put(1651,-436){\makebox(0,0)[lb]{\smash{\SetFigFont{8}{9.6}{\rmdefault}{\mddefault}{\updefault}1}}}
\put(2326,-436){\makebox(0,0)[lb]{\smash{\SetFigFont{8}{9.6}{\rmdefault}{\mddefault}{\updefault}1}}}
\put(1501,-961){\makebox(0,0)[lb]{\smash{\SetFigFont{8}{9.6}{\rmdefault}{\mddefault}{\updefault}a}}}
\put(1951,-961){\makebox(0,0)[lb]{\smash{\SetFigFont{8}{9.6}{\rmdefault}{\mddefault}{\updefault}b}}}
\end{picture}} 
\parbox{2cm}{$\; = \; \frac{A^{-3}}{[2]}$}
\parbox{2cm}{\setlength{\unitlength}{0.00050000in}%
\begingroup\makeatletter\ifx\SetFigFont\undefined%
\gdef\SetFigFont#1#2#3#4#5{%
  \reset@font\fontsize{#1}{#2pt}%
  \fontfamily{#3}\fontseries{#4}\fontshape{#5}%
  \selectfont}%
\fi\endgroup%
\begin{picture}(1074,1110)(1564,-1063)
\thicklines
\put(1576,-736){\line( 1, 0){1050}}
\multiput(1876,-736)(10.22727,10.22727){23}{\makebox(11.1111,16.6667){\SetFigFont{7}{8.4}{\rmdefault}{\mddefault}{\updefault}.}}
\multiput(2101,-511)(10.22727,-10.22727){23}{\makebox(11.1111,16.6667){\SetFigFont{7}{8.4}{\rmdefault}{\mddefault}{\updefault}.}}
\put(2101,-511){\line( 0, 1){300}}
\multiput(2101,-211)(10.22727,10.22727){23}{\makebox(11.1111,16.6667){\SetFigFont{7}{8.4}{\rmdefault}{\mddefault}{\updefault}.}}
\multiput(1876, 14)(10.22727,-10.22727){23}{\makebox(11.1111,16.6667){\SetFigFont{7}{8.4}{\rmdefault}{\mddefault}{\updefault}.}}
\put(1576,-1036){\makebox(0,0)[lb]{\smash{\SetFigFont{7}{8.4}{\rmdefault}{\mddefault}{\updefault}a}}}
\put(2401,-61){\makebox(0,0)[lb]{\smash{\SetFigFont{7}{8.4}{\rmdefault}{\mddefault}{\updefault}1}}}
\put(1726,-61){\makebox(0,0)[lb]{\smash{\SetFigFont{7}{8.4}{\rmdefault}{\mddefault}{\updefault}1}}}
\put(1876,-586){\makebox(0,0)[lb]{\smash{\SetFigFont{7}{8.4}{\rmdefault}{\mddefault}{\updefault}1}}}
\put(2326,-586){\makebox(0,0)[lb]{\smash{\SetFigFont{7}{8.4}{\rmdefault}{\mddefault}{\updefault}1}}}
\put(2026,-1036){\makebox(0,0)[lb]{\smash{\SetFigFont{7}{8.4}{\rmdefault}{\mddefault}{\updefault}b}}}
\put(2476,-1036){\makebox(0,0)[lb]{\smash{\SetFigFont{7}{8.4}{\rmdefault}{\mddefault}{\updefault}a}}}
\put(2176,-361){\makebox(0,0)[lb]{\smash{\SetFigFont{7}{8.4}{\rmdefault}{\mddefault}{\updefault}0}}}
\end{picture}}
\parbox{1cm}{$\; + \; A$}
\parbox{2cm}{\setlength{\unitlength}{0.00050000in}%
\begingroup\makeatletter\ifx\SetFigFont\undefined%
\gdef\SetFigFont#1#2#3#4#5{%
  \reset@font\fontsize{#1}{#2pt}%
  \fontfamily{#3}\fontseries{#4}\fontshape{#5}%
  \selectfont}%
\fi\endgroup%
\begin{picture}(1074,1110)(1564,-1063)
\thicklines
\put(1576,-736){\line( 1, 0){1050}}
\multiput(1876,-736)(10.22727,10.22727){23}{\makebox(11.1111,16.6667){\SetFigFont{7}{8.4}{\rmdefault}{\mddefault}{\updefault}.}}
\multiput(2101,-511)(10.22727,-10.22727){23}{\makebox(11.1111,16.6667){\SetFigFont{7}{8.4}{\rmdefault}{\mddefault}{\updefault}.}}
\put(2101,-511){\line( 0, 1){300}}
\multiput(2101,-211)(10.22727,10.22727){23}{\makebox(11.1111,16.6667){\SetFigFont{7}{8.4}{\rmdefault}{\mddefault}{\updefault}.}}
\multiput(1876, 14)(10.22727,-10.22727){23}{\makebox(11.1111,16.6667){\SetFigFont{7}{8.4}{\rmdefault}{\mddefault}{\updefault}.}}
\put(1576,-1036){\makebox(0,0)[lb]{\smash{\SetFigFont{7}{8.4}{\rmdefault}{\mddefault}{\updefault}a}}}
\put(2401,-61){\makebox(0,0)[lb]{\smash{\SetFigFont{7}{8.4}{\rmdefault}{\mddefault}{\updefault}1}}}
\put(1726,-61){\makebox(0,0)[lb]{\smash{\SetFigFont{7}{8.4}{\rmdefault}{\mddefault}{\updefault}1}}}
\put(1876,-586){\makebox(0,0)[lb]{\smash{\SetFigFont{7}{8.4}{\rmdefault}{\mddefault}{\updefault}1}}}
\put(2326,-586){\makebox(0,0)[lb]{\smash{\SetFigFont{7}{8.4}{\rmdefault}{\mddefault}{\updefault}1}}}
\put(2176,-361){\makebox(0,0)[lb]{\smash{\SetFigFont{7}{8.4}{\rmdefault}{\mddefault}{\updefault}2}}}
\put(2026,-1036){\makebox(0,0)[lb]{\smash{\SetFigFont{7}{8.4}{\rmdefault}{\mddefault}{\updefault}b}}}
\put(2476,-1036){\makebox(0,0)[lb]{\smash{\SetFigFont{7}{8.4}{\rmdefault}{\mddefault}{\updefault}a}}}
\end{picture}}
\parbox{1cm}{$\; = \;$}
\end{center}

\begin{center}
\parbox{1cm}{$\; = \;$}
\parbox{4cm}{$\; \; \frac{A^{-3}}{[2]}
 \frac{\left < 
                \begin{array}{ccc}
                 0 & 1 & 1 \\
                 b & a & a
                 \end{array} \right >}{<0,a,a>}$}
\parbox{2cm}{\setlength{\unitlength}{0.00050000in}%
\begingroup\makeatletter\ifx\SetFigFont\undefined%
\gdef\SetFigFont#1#2#3#4#5{%
  \reset@font\fontsize{#1}{#2pt}%
  \fontfamily{#3}\fontseries{#4}\fontshape{#5}%
  \selectfont}%
\fi\endgroup%
\begin{picture}(762,855)(1726,-1054)
\thicklines
\put(2101,-811){\makebox(11.1111,16.6667){\SetFigFont{10}{12}{\rmdefault}{\mddefault}{\updefault}.}}
\put(2101,-811){\line( 0, 1){450}}
\multiput(2101,-361)(-11.84211,7.89474){20}{\makebox(11.1111,16.6667){\SetFigFont{7}{8.4}{\rmdefault}{\mddefault}{\updefault}.}}
\put(1801,-811){\line( 1, 0){675}}
\multiput(2101,-361)(11.84211,7.89474){20}{\makebox(11.1111,16.6667){\SetFigFont{7}{8.4}{\rmdefault}{\mddefault}{\updefault}.}}
\put(1876,-1036){\makebox(0,0)[lb]{\smash{\SetFigFont{7}{8.4}{\rmdefault}{\mddefault}{\updefault}a}}}
\put(1726,-361){\makebox(0,0)[lb]{\smash{\SetFigFont{7}{8.4}{\rmdefault}{\mddefault}{\updefault}1}}}
\put(2401,-361){\makebox(0,0)[lb]{\smash{\SetFigFont{7}{8.4}{\rmdefault}{\mddefault}{\updefault}1}}}
\put(2251,-1036){\makebox(0,0)[lb]{\smash{\SetFigFont{7}{8.4}{\rmdefault}{\mddefault}{\updefault}a}}}
\put(2176,-661){\makebox(0,0)[lb]{\smash{\SetFigFont{7}{8.4}{\rmdefault}{\mddefault}{\updefault}0}}}
\end{picture}}
\parbox{1cm}{$\; + \;$}
\parbox{4cm}{$\;  \; A
 \frac{\left < 
                \begin{array}{ccc}
                 2 & 1 & 1 \\
                 b & a & a
                 \end{array} \right >}{<2,a,a>}$}
\parbox{2cm}{\setlength{\unitlength}{0.00050000in}%
\begingroup\makeatletter\ifx\SetFigFont\undefined%
\gdef\SetFigFont#1#2#3#4#5{%
  \reset@font\fontsize{#1}{#2pt}%
  \fontfamily{#3}\fontseries{#4}\fontshape{#5}%
  \selectfont}%
\fi\endgroup%
\begin{picture}(762,855)(1726,-1054)
\thicklines
\put(2101,-811){\makebox(11.1111,16.6667){\SetFigFont{10}{12}{\rmdefault}{\mddefault}{\updefault}.}}
\put(2101,-811){\line( 0, 1){450}}
\multiput(2101,-361)(-11.84211,7.89474){20}{\makebox(11.1111,16.6667){\SetFigFont{7}{8.4}{\rmdefault}{\mddefault}{\updefault}.}}
\put(1801,-811){\line( 1, 0){675}}
\multiput(2101,-361)(11.84211,7.89474){20}{\makebox(11.1111,16.6667){\SetFigFont{7}{8.4}{\rmdefault}{\mddefault}{\updefault}.}}
\put(1876,-1036){\makebox(0,0)[lb]{\smash{\SetFigFont{7}{8.4}{\rmdefault}{\mddefault}{\updefault}a}}}
\put(2176,-661){\makebox(0,0)[lb]{\smash{\SetFigFont{7}{8.4}{\rmdefault}{\mddefault}{\updefault}2}}}
\put(1726,-361){\makebox(0,0)[lb]{\smash{\SetFigFont{7}{8.4}{\rmdefault}{\mddefault}{\updefault}1}}}
\put(2401,-361){\makebox(0,0)[lb]{\smash{\SetFigFont{7}{8.4}{\rmdefault}{\mddefault}{\updefault}1}}}
\put(2251,-1036){\makebox(0,0)[lb]{\smash{\SetFigFont{7}{8.4}{\rmdefault}{\mddefault}{\updefault}a}}}
\end{picture}}
\end{center}
\vspace{0.5cm} 

\noindent We have to perform a fusing, in order to arrive to the standard
basis of $V$. Using the computations of 6j-symbols appearing in this
particular fusing we obtain that:

\begin{center} 
\parbox{2cm}{\setlength{\unitlength}{0.00050000in}%
\begingroup\makeatletter\ifx\SetFigFont\undefined%
\gdef\SetFigFont#1#2#3#4#5{%
  \reset@font\fontsize{#1}{#2pt}%
  \fontfamily{#3}\fontseries{#4}\fontshape{#5}%
  \selectfont}%
\fi\endgroup%
\begin{picture}(762,855)(1726,-1054)
\thicklines
\put(2101,-811){\makebox(11.1111,16.6667){\SetFigFont{10}{12}{\rmdefault}{\mddefault}{\updefault}.}}
\put(2101,-811){\line( 0, 1){450}}
\multiput(2101,-361)(-11.84211,7.89474){20}{\makebox(11.1111,16.6667){\SetFigFont{7}{8.4}{\rmdefault}{\mddefault}{\updefault}.}}
\put(1801,-811){\line( 1, 0){675}}
\multiput(2101,-361)(11.84211,7.89474){20}{\makebox(11.1111,16.6667){\SetFigFont{7}{8.4}{\rmdefault}{\mddefault}{\updefault}.}}
\put(1876,-1036){\makebox(0,0)[lb]{\smash{\SetFigFont{7}{8.4}{\rmdefault}{\mddefault}{\updefault}a}}}
\put(1726,-361){\makebox(0,0)[lb]{\smash{\SetFigFont{7}{8.4}{\rmdefault}{\mddefault}{\updefault}1}}}
\put(2401,-361){\makebox(0,0)[lb]{\smash{\SetFigFont{7}{8.4}{\rmdefault}{\mddefault}{\updefault}1}}}
\put(2251,-1036){\makebox(0,0)[lb]{\smash{\SetFigFont{7}{8.4}{\rmdefault}{\mddefault}{\updefault}a}}}
\put(2176,-661){\makebox(0,0)[lb]{\smash{\SetFigFont{7}{8.4}{\rmdefault}{\mddefault}{\updefault}0}}}
\end{picture}} 
\parbox{2cm}{$\; = \; -\frac{[a]}{[a+1]}$}
\parbox{2cm}{\setlength{\unitlength}{0.00050000in}%
\begingroup\makeatletter\ifx\SetFigFont\undefined%
\gdef\SetFigFont#1#2#3#4#5{%
  \reset@font\fontsize{#1}{#2pt}%
  \fontfamily{#3}\fontseries{#4}\fontshape{#5}%
  \selectfont}%
\fi\endgroup%
\begin{picture}(1074,714)(1564,-1063)
\thicklines
\put(1576,-736){\line( 1, 0){1050}}
\put(1876,-736){\line( 0, 1){375}}
\put(2251,-736){\line( 0, 1){375}}
\put(1651,-1036){\makebox(0,0)[lb]{\smash{\SetFigFont{7}{8.4}{\rmdefault}{\mddefault}{\updefault}a}}}
\put(1726,-511){\makebox(0,0)[lb]{\smash{\SetFigFont{7}{8.4}{\rmdefault}{\mddefault}{\updefault}1}}}
\put(2326,-511){\makebox(0,0)[lb]{\smash{\SetFigFont{7}{8.4}{\rmdefault}{\mddefault}{\updefault}1}}}
\put(1951,-1036){\makebox(0,0)[lb]{\smash{\SetFigFont{7}{8.4}{\rmdefault}{\mddefault}{\updefault}a-1}}}
\put(2476,-1036){\makebox(0,0)[lb]{\smash{\SetFigFont{7}{8.4}{\rmdefault}{\mddefault}{\updefault}a}}}
\end{picture}}
\parbox{1cm}{$\; + \; $}
\parbox{2cm}{\setlength{\unitlength}{0.00050000in}%
\begingroup\makeatletter\ifx\SetFigFont\undefined%
\gdef\SetFigFont#1#2#3#4#5{%
  \reset@font\fontsize{#1}{#2pt}%
  \fontfamily{#3}\fontseries{#4}\fontshape{#5}%
  \selectfont}%
\fi\endgroup%
\begin{picture}(1074,714)(1564,-1063)
\thicklines
\put(1576,-736){\line( 1, 0){1050}}
\put(1876,-736){\line( 0, 1){375}}
\put(2251,-736){\line( 0, 1){375}}
\put(1651,-1036){\makebox(0,0)[lb]{\smash{\SetFigFont{7}{8.4}{\rmdefault}{\mddefault}{\updefault}a}}}
\put(1951,-1036){\makebox(0,0)[lb]{\smash{\SetFigFont{7}{8.4}{\rmdefault}{\mddefault}{\updefault}a+1}}}
\put(1726,-511){\makebox(0,0)[lb]{\smash{\SetFigFont{7}{8.4}{\rmdefault}{\mddefault}{\updefault}1}}}
\put(2326,-511){\makebox(0,0)[lb]{\smash{\SetFigFont{7}{8.4}{\rmdefault}{\mddefault}{\updefault}1}}}
\put(2401,-1036){\makebox(0,0)[lb]{\smash{\SetFigFont{7}{8.4}{\rmdefault}{\mddefault}{\updefault}a}}}
\end{picture}}
\end{center}

\begin{center} 
\parbox{2cm}{\setlength{\unitlength}{0.00050000in}%
\begingroup\makeatletter\ifx\SetFigFont\undefined%
\gdef\SetFigFont#1#2#3#4#5{%
  \reset@font\fontsize{#1}{#2pt}%
  \fontfamily{#3}\fontseries{#4}\fontshape{#5}%
  \selectfont}%
\fi\endgroup%
\begin{picture}(762,855)(1726,-1054)
\thicklines
\put(2101,-811){\makebox(11.1111,16.6667){\SetFigFont{10}{12}{\rmdefault}{\mddefault}{\updefault}.}}
\put(2101,-811){\line( 0, 1){450}}
\multiput(2101,-361)(-11.84211,7.89474){20}{\makebox(11.1111,16.6667){\SetFigFont{7}{8.4}{\rmdefault}{\mddefault}{\updefault}.}}
\put(1801,-811){\line( 1, 0){675}}
\multiput(2101,-361)(11.84211,7.89474){20}{\makebox(11.1111,16.6667){\SetFigFont{7}{8.4}{\rmdefault}{\mddefault}{\updefault}.}}
\put(1876,-1036){\makebox(0,0)[lb]{\smash{\SetFigFont{7}{8.4}{\rmdefault}{\mddefault}{\updefault}a}}}
\put(2176,-661){\makebox(0,0)[lb]{\smash{\SetFigFont{7}{8.4}{\rmdefault}{\mddefault}{\updefault}2}}}
\put(1726,-361){\makebox(0,0)[lb]{\smash{\SetFigFont{7}{8.4}{\rmdefault}{\mddefault}{\updefault}1}}}
\put(2401,-361){\makebox(0,0)[lb]{\smash{\SetFigFont{7}{8.4}{\rmdefault}{\mddefault}{\updefault}1}}}
\put(2251,-1036){\makebox(0,0)[lb]{\smash{\SetFigFont{7}{8.4}{\rmdefault}{\mddefault}{\updefault}a}}}
\end{picture}} 
\parbox{2cm}{$\; = \; \frac{[a+2]}{[2][a+1]}$}
\parbox{2cm}{\setlength{\unitlength}{0.00050000in}%
\begingroup\makeatletter\ifx\SetFigFont\undefined%
\gdef\SetFigFont#1#2#3#4#5{%
  \reset@font\fontsize{#1}{#2pt}%
  \fontfamily{#3}\fontseries{#4}\fontshape{#5}%
  \selectfont}%
\fi\endgroup%
\begin{picture}(1074,714)(1564,-1063)
\thicklines
\put(1576,-736){\line( 1, 0){1050}}
\put(1876,-736){\line( 0, 1){375}}
\put(2251,-736){\line( 0, 1){375}}
\put(1651,-1036){\makebox(0,0)[lb]{\smash{\SetFigFont{7}{8.4}{\rmdefault}{\mddefault}{\updefault}a}}}
\put(1726,-511){\makebox(0,0)[lb]{\smash{\SetFigFont{7}{8.4}{\rmdefault}{\mddefault}{\updefault}1}}}
\put(2326,-511){\makebox(0,0)[lb]{\smash{\SetFigFont{7}{8.4}{\rmdefault}{\mddefault}{\updefault}1}}}
\put(1951,-1036){\makebox(0,0)[lb]{\smash{\SetFigFont{7}{8.4}{\rmdefault}{\mddefault}{\updefault}a-1}}}
\put(2476,-1036){\makebox(0,0)[lb]{\smash{\SetFigFont{7}{8.4}{\rmdefault}{\mddefault}{\updefault}a}}}
\end{picture}}
\parbox{1cm}{$\; + \; \frac{1}{[2]}$}
\parbox{2cm}{\setlength{\unitlength}{0.00050000in}%
\begingroup\makeatletter\ifx\SetFigFont\undefined%
\gdef\SetFigFont#1#2#3#4#5{%
  \reset@font\fontsize{#1}{#2pt}%
  \fontfamily{#3}\fontseries{#4}\fontshape{#5}%
  \selectfont}%
\fi\endgroup%
\begin{picture}(1074,714)(1564,-1063)
\thicklines
\put(1576,-736){\line( 1, 0){1050}}
\put(1876,-736){\line( 0, 1){375}}
\put(2251,-736){\line( 0, 1){375}}
\put(1651,-1036){\makebox(0,0)[lb]{\smash{\SetFigFont{7}{8.4}{\rmdefault}{\mddefault}{\updefault}a}}}
\put(1951,-1036){\makebox(0,0)[lb]{\smash{\SetFigFont{7}{8.4}{\rmdefault}{\mddefault}{\updefault}a+1}}}
\put(1726,-511){\makebox(0,0)[lb]{\smash{\SetFigFont{7}{8.4}{\rmdefault}{\mddefault}{\updefault}1}}}
\put(2326,-511){\makebox(0,0)[lb]{\smash{\SetFigFont{7}{8.4}{\rmdefault}{\mddefault}{\updefault}1}}}
\put(2401,-1036){\makebox(0,0)[lb]{\smash{\SetFigFont{7}{8.4}{\rmdefault}{\mddefault}{\updefault}a}}}
\end{picture}}
\end{center}

\noindent The convention is that diagrams whose labels form non-admissible
triples are vanishing. Explicit computations now yield our claim.
$\Box$

\newtheorem{lem4}[lem1]{Lemma}
\begin{lem4}\label{lem4}
The following numerical identities are satisfied:
\[ [a]+[a+2] = [2][a+1], 
\; \frac{[a+1](1+A^{-4})}{[a]+[a+2]}= \frac{1}{A^2},\]
\[\frac{[a+2]A^{-4}-[a]}{[2][a+1]}= \frac{A^{-4}-1}{1-A^{4+4a}}, 
\; \frac{[a]A^{-4}-[a+2]}{[2][a+1]}= \frac{A^{-4}-1}{1-A^{-4-4a}}. \]
\end{lem4}
The proof is a mere computation.$\Box$

\noindent Set now $q= A^{-4}$. Then lemma 3 can be reformulated as 

\begin{center} 
\parbox{1cm}{$\rho(g_i)$}
\parbox{2cm}{\setlength{\unitlength}{0.00050000in}%
\begingroup\makeatletter\ifx\SetFigFont\undefined%
\gdef\SetFigFont#1#2#3#4#5{%
  \reset@font\fontsize{#1}{#2pt}%
  \fontfamily{#3}\fontseries{#4}\fontshape{#5}%
  \selectfont}%
\fi\endgroup%
\begin{picture}(1074,714)(1564,-1063)
\thicklines
\put(1576,-736){\line( 1, 0){1050}}
\put(1876,-736){\line( 0, 1){375}}
\put(2251,-736){\line( 0, 1){375}}
\put(1651,-1036){\makebox(0,0)[lb]{\smash{\SetFigFont{7}{8.4}{\rmdefault}{\mddefault}{\updefault}a}}}
\put(1726,-511){\makebox(0,0)[lb]{\smash{\SetFigFont{7}{8.4}{\rmdefault}{\mddefault}{\updefault}1}}}
\put(2326,-511){\makebox(0,0)[lb]{\smash{\SetFigFont{7}{8.4}{\rmdefault}{\mddefault}{\updefault}1}}}
\put(1951,-1036){\makebox(0,0)[lb]{\smash{\SetFigFont{7}{8.4}{\rmdefault}{\mddefault}{\updefault}a-1}}}
\put(2476,-1036){\makebox(0,0)[lb]{\smash{\SetFigFont{7}{8.4}{\rmdefault}{\mddefault}{\updefault}a}}}
\end{picture}} 
\parbox{3cm}{$\; = \; - A \frac{q-1}{1-q^{a+1}}$}
\parbox{2cm}{\setlength{\unitlength}{0.00050000in}%
\begingroup\makeatletter\ifx\SetFigFont\undefined%
\gdef\SetFigFont#1#2#3#4#5{%
  \reset@font\fontsize{#1}{#2pt}%
  \fontfamily{#3}\fontseries{#4}\fontshape{#5}%
  \selectfont}%
\fi\endgroup%
\begin{picture}(1074,714)(1564,-1063)
\thicklines
\put(1576,-736){\line( 1, 0){1050}}
\put(1876,-736){\line( 0, 1){375}}
\put(2251,-736){\line( 0, 1){375}}
\put(1651,-1036){\makebox(0,0)[lb]{\smash{\SetFigFont{7}{8.4}{\rmdefault}{\mddefault}{\updefault}a}}}
\put(1726,-511){\makebox(0,0)[lb]{\smash{\SetFigFont{7}{8.4}{\rmdefault}{\mddefault}{\updefault}1}}}
\put(2326,-511){\makebox(0,0)[lb]{\smash{\SetFigFont{7}{8.4}{\rmdefault}{\mddefault}{\updefault}1}}}
\put(1951,-1036){\makebox(0,0)[lb]{\smash{\SetFigFont{7}{8.4}{\rmdefault}{\mddefault}{\updefault}a-1}}}
\put(2476,-1036){\makebox(0,0)[lb]{\smash{\SetFigFont{7}{8.4}{\rmdefault}{\mddefault}{\updefault}a}}}
\end{picture}}
\parbox{2cm}{$\; +\; \frac{1}{A}  \; $}
\parbox{2cm}{\setlength{\unitlength}{0.00050000in}%
\begingroup\makeatletter\ifx\SetFigFont\undefined%
\gdef\SetFigFont#1#2#3#4#5{%
  \reset@font\fontsize{#1}{#2pt}%
  \fontfamily{#3}\fontseries{#4}\fontshape{#5}%
  \selectfont}%
\fi\endgroup%
\begin{picture}(1074,714)(1564,-1063)
\thicklines
\put(1576,-736){\line( 1, 0){1050}}
\put(1876,-736){\line( 0, 1){375}}
\put(2251,-736){\line( 0, 1){375}}
\put(1651,-1036){\makebox(0,0)[lb]{\smash{\SetFigFont{7}{8.4}{\rmdefault}{\mddefault}{\updefault}a}}}
\put(1951,-1036){\makebox(0,0)[lb]{\smash{\SetFigFont{7}{8.4}{\rmdefault}{\mddefault}{\updefault}a+1}}}
\put(1726,-511){\makebox(0,0)[lb]{\smash{\SetFigFont{7}{8.4}{\rmdefault}{\mddefault}{\updefault}1}}}
\put(2326,-511){\makebox(0,0)[lb]{\smash{\SetFigFont{7}{8.4}{\rmdefault}{\mddefault}{\updefault}1}}}
\put(2401,-1036){\makebox(0,0)[lb]{\smash{\SetFigFont{7}{8.4}{\rmdefault}{\mddefault}{\updefault}a}}}
\end{picture}}
\end{center}

\begin{center} 
\parbox{1cm}{$\rho(g_i)$}
\parbox{2cm}{\setlength{\unitlength}{0.00050000in}%
\begingroup\makeatletter\ifx\SetFigFont\undefined%
\gdef\SetFigFont#1#2#3#4#5{%
  \reset@font\fontsize{#1}{#2pt}%
  \fontfamily{#3}\fontseries{#4}\fontshape{#5}%
  \selectfont}%
\fi\endgroup%
\begin{picture}(1074,714)(1564,-1063)
\thicklines
\put(1576,-736){\line( 1, 0){1050}}
\put(1876,-736){\line( 0, 1){375}}
\put(2251,-736){\line( 0, 1){375}}
\put(1651,-1036){\makebox(0,0)[lb]{\smash{\SetFigFont{7}{8.4}{\rmdefault}{\mddefault}{\updefault}a}}}
\put(1951,-1036){\makebox(0,0)[lb]{\smash{\SetFigFont{7}{8.4}{\rmdefault}{\mddefault}{\updefault}a+1}}}
\put(1726,-511){\makebox(0,0)[lb]{\smash{\SetFigFont{7}{8.4}{\rmdefault}{\mddefault}{\updefault}1}}}
\put(2326,-511){\makebox(0,0)[lb]{\smash{\SetFigFont{7}{8.4}{\rmdefault}{\mddefault}{\updefault}1}}}
\put(2401,-1036){\makebox(0,0)[lb]{\smash{\SetFigFont{7}{8.4}{\rmdefault}{\mddefault}{\updefault}a}}}
\end{picture}} 
\parbox{2cm}{$\; = \; 
\nu^2A^3 \; $}
\parbox{2cm}{\setlength{\unitlength}{0.00050000in}%
\begingroup\makeatletter\ifx\SetFigFont\undefined%
\gdef\SetFigFont#1#2#3#4#5{%
  \reset@font\fontsize{#1}{#2pt}%
  \fontfamily{#3}\fontseries{#4}\fontshape{#5}%
  \selectfont}%
\fi\endgroup%
\begin{picture}(1074,714)(1564,-1063)
\thicklines
\put(1576,-736){\line( 1, 0){1050}}
\put(1876,-736){\line( 0, 1){375}}
\put(2251,-736){\line( 0, 1){375}}
\put(1651,-1036){\makebox(0,0)[lb]{\smash{\SetFigFont{7}{8.4}{\rmdefault}{\mddefault}{\updefault}a}}}
\put(1726,-511){\makebox(0,0)[lb]{\smash{\SetFigFont{7}{8.4}{\rmdefault}{\mddefault}{\updefault}1}}}
\put(2326,-511){\makebox(0,0)[lb]{\smash{\SetFigFont{7}{8.4}{\rmdefault}{\mddefault}{\updefault}1}}}
\put(1951,-1036){\makebox(0,0)[lb]{\smash{\SetFigFont{7}{8.4}{\rmdefault}{\mddefault}{\updefault}a-1}}}
\put(2476,-1036){\makebox(0,0)[lb]{\smash{\SetFigFont{7}{8.4}{\rmdefault}{\mddefault}{\updefault}a}}}
\end{picture}}
\parbox{3cm}{$\; - \; A  \frac{q-1}{1-q^{-a-1}} \; $}
\parbox{2cm}{\setlength{\unitlength}{0.00050000in}%
\begingroup\makeatletter\ifx\SetFigFont\undefined%
\gdef\SetFigFont#1#2#3#4#5{%
  \reset@font\fontsize{#1}{#2pt}%
  \fontfamily{#3}\fontseries{#4}\fontshape{#5}%
  \selectfont}%
\fi\endgroup%
\begin{picture}(1074,714)(1564,-1063)
\thicklines
\put(1576,-736){\line( 1, 0){1050}}
\put(1876,-736){\line( 0, 1){375}}
\put(2251,-736){\line( 0, 1){375}}
\put(1651,-1036){\makebox(0,0)[lb]{\smash{\SetFigFont{7}{8.4}{\rmdefault}{\mddefault}{\updefault}a}}}
\put(1951,-1036){\makebox(0,0)[lb]{\smash{\SetFigFont{7}{8.4}{\rmdefault}{\mddefault}{\updefault}a+1}}}
\put(1726,-511){\makebox(0,0)[lb]{\smash{\SetFigFont{7}{8.4}{\rmdefault}{\mddefault}{\updefault}1}}}
\put(2326,-511){\makebox(0,0)[lb]{\smash{\SetFigFont{7}{8.4}{\rmdefault}{\mddefault}{\updefault}1}}}
\put(2401,-1036){\makebox(0,0)[lb]{\smash{\SetFigFont{7}{8.4}{\rmdefault}{\mddefault}{\updefault}a}}}
\end{picture}}
\end{center}

\noindent where $\nu^2=\frac{(1-q^{a+2})(1-q^{a+1})(1-q^{-a})(1-q^{-a-1})}
{(1+q)^2}$. 

We specify now to the cases $n=3, \, m=1$, and $n=4, m=2$. 
Then the vector  space $V(n,m)$ has dimension 2 and respectively 3. 
The Hecke algebra relations follow immediately from the formulas 
above. The additional relation
defining $A_{\beta,n}$ is  verified by a direct computation. 

Let us check out the case $n=3$. The two vectors which span 
$V(3,1)$ have the labels $(0,p_1,p_2,m)\in \{ (0,1,0,1),
(0,1,2,1)\}$. Then the previous formulas read:
\[ \tilde{\rho}(g_1) = 
\left ( \begin{array}{cc}
-1 & 0 \\
0 & A^{-4}
\end{array}
\right ), \]
\[ \tilde{\rho}(g_2) = 
\left ( \begin{array}{cc}
\frac{1}{A^4(1+A^4)} & -\frac{A^4+A^{-4}+1}{A^2(A^4+A^{-4}+2)} \\
-A^{-2} & -\frac{1}{1+A^{-4}}
\end{array}
\right ),  \]
where $\tilde{\rho}=(-A^{-1})\rho$. This implies that 
\[ \tilde{\rho}(g_1g_2) = 
\left ( \begin{array}{cc}
-\frac{1}{A^4(1+A^4)} & \frac{A^4+A^{-4}+1}{A^2(A^4+A^{-4}+2)} \\
-A^{-6} & -\frac{1}{A^4(1+A^{-4})}
\end{array}
\right ),  \; 
 \tilde{\rho}(g_2g_1) = 
\left ( \begin{array}{cc}
-\frac{1}{A^4(1+A^4)} & -\frac{A^4+A^{-4}+1}{A^6(A^4+A^{-4}+2)} \\
A^{-2} & -\frac{1}{A^4(1+A^{-4})}
\end{array}
\right ),  \]
and 
\[ \tilde{\rho}(g_1g_2g_1) = 
\left ( \begin{array}{cc}
\frac{1}{A^4(1+A^4)} & \frac{A^4+A^{-4}+1}{A^6(A^4+A^{-4}+2)} \\
A^{-6} & -\frac{1}{A^8(1+A^{-4})}
\end{array}
\right ).  \]
Therefore the relation  
$1+\tilde{\rho}(g_1)+\tilde{\rho}(g_2)+\tilde{\rho}(g_1g_2)+
\tilde{\rho}(g_2g_1)+\tilde{\rho}(g_1g_2g_1)=0$
holds, which proves the first part of  proposition \ref{pr}.

Let us give the explicit matrices for $n=4, \; m=2$ and an 
ad-hoc proof of the irreducibility. 
The space $V(4,2)$ is spanned by the vectors of labels 
$(0,p_1,p_2,p_3,m)\in \{ (0,1,0,1,2), (0,1,2,1,2), (0,1,2,3,2)\}$.
Then within this basis we have from above:
\[ \tilde{\rho}(g_1) = 
\left ( \begin{array}{ccc}
A^{-4} & 0 & 0 \\
0 & -1 &  0 \\
0 & 0 & -1 
\end{array}
\right ),  \]
\[ \tilde{\rho}(g_2) = 
\left ( \begin{array}{ccc}
-\frac{A}{1+A^{-4}} & -A^{-2} & 0 \\
-\frac{(1-A^{12})(1-A^4)^3}{A^{34}} & \frac{1}{1+A^{8}} &  0 \\
0 & 0 & -1 
\end{array}
\right ),  \]
\[ \tilde{\rho}(g_3) = 
\left ( \begin{array}{ccc}
-1 & 0 & 0 \\
0 & -\frac{1}{1+A^{-4}+A^{-8}} &  -A^{-2} \\
0 & -\frac{(1+A^{8})^2(A^{4}-1)^2(A^{12}-1)^2}{A^{66}} & \frac{A^{-12} }{1+A^{-4}+A^{-8}}
\end{array}
\right ).   \]
Assume that the $B_4$ representation $\tilde{\rho}$ is not irreducible. 
Then the $A_{\beta,4}$-module $V(4,2)$ is completely reducible and
there is at least one simple factor of 
dimension 1. Equivalently,  $V(4,2)$  contains  a 1-dimensional
$B_4$-invariant subspace say ${\bf C}w$, for some 
non-zero vector $w$. There exist then the scalars 
$\lambda_i$ such that  $\tilde{\rho}(g_i)w =
 \lambda_i w$. The group relations imply
 $\lambda_1=\lambda_2=\lambda_3=\lambda$ (since the matrices above are
 non-singular) and from the relations in $A_{\beta,4}$ we derive 
$\lambda=-1$. However the condition $\tilde{\rho}(g_1)w+w=0$ 
yields $w_1=0$, further $\tilde{\rho}(g_2)w+w=0$ adds the 
constraint $w_2=0$ and the last identity shows that $w$ vanishes. 
This proves that $V(4,2)$ is irreducible. 
In particular $\tilde{\rho}$ factors throughout 
the projection $A_{\beta,4}\longrightarrow M_3({\bf C})$.  
$\Box$
\newtheorem{rem}[lem1]{Remarks}
\begin{rem}
 We could show from the very beginning of the proof that the
  representation of $B_n$ factors through the Hecke algebra $H_n(q)$ 
with $q=A^{-4}$. 
Let assume we are interested in the action of  $g_j$.  
Observe  that the  vectors

\begin{center} 
\parbox{2cm}{\setlength{\unitlength}{0.00050000in}%
\begingroup\makeatletter\ifx\SetFigFont\undefined%
\gdef\SetFigFont#1#2#3#4#5{%
  \reset@font\fontsize{#1}{#2pt}%
  \fontfamily{#3}\fontseries{#4}\fontshape{#5}%
  \selectfont}%
\fi\endgroup%
\begin{picture}(762,855)(1726,-1054)
\thicklines
\put(2101,-811){\makebox(11.1111,16.6667){\SetFigFont{10}{12}{\rmdefault}{\mddefault}{\updefault}.}}
\put(2101,-811){\line( 0, 1){450}}
\multiput(2101,-361)(-11.84211,7.89474){20}{\makebox(11.1111,16.6667){\SetFigFont{7}{8.4}{\rmdefault}{\mddefault}{\updefault}.}}
\put(1801,-811){\line( 1, 0){675}}
\multiput(2101,-361)(11.84211,7.89474){20}{\makebox(11.1111,16.6667){\SetFigFont{7}{8.4}{\rmdefault}{\mddefault}{\updefault}.}}
\put(1876,-1036){\makebox(0,0)[lb]{\smash{\SetFigFont{7}{8.4}{\rmdefault}{\mddefault}{\updefault}a}}}
\put(1726,-361){\makebox(0,0)[lb]{\smash{\SetFigFont{7}{8.4}{\rmdefault}{\mddefault}{\updefault}1}}}
\put(2401,-361){\makebox(0,0)[lb]{\smash{\SetFigFont{7}{8.4}{\rmdefault}{\mddefault}{\updefault}1}}}
\put(2251,-1036){\makebox(0,0)[lb]{\smash{\SetFigFont{7}{8.4}{\rmdefault}{\mddefault}{\updefault}a}}}
\put(2176,-661){\makebox(0,0)[lb]{\smash{\SetFigFont{7}{8.4}{\rmdefault}{\mddefault}{\updefault}0}}}
\end{picture}}
\parbox{0.4cm}{,}
\parbox{2cm}{\setlength{\unitlength}{0.00050000in}%
\begingroup\makeatletter\ifx\SetFigFont\undefined%
\gdef\SetFigFont#1#2#3#4#5{%
  \reset@font\fontsize{#1}{#2pt}%
  \fontfamily{#3}\fontseries{#4}\fontshape{#5}%
  \selectfont}%
\fi\endgroup%
\begin{picture}(762,855)(1726,-1054)
\thicklines
\put(2101,-811){\makebox(11.1111,16.6667){\SetFigFont{10}{12}{\rmdefault}{\mddefault}{\updefault}.}}
\put(2101,-811){\line( 0, 1){450}}
\multiput(2101,-361)(-11.84211,7.89474){20}{\makebox(11.1111,16.6667){\SetFigFont{7}{8.4}{\rmdefault}{\mddefault}{\updefault}.}}
\put(1801,-811){\line( 1, 0){675}}
\multiput(2101,-361)(11.84211,7.89474){20}{\makebox(11.1111,16.6667){\SetFigFont{7}{8.4}{\rmdefault}{\mddefault}{\updefault}.}}
\put(1876,-1036){\makebox(0,0)[lb]{\smash{\SetFigFont{7}{8.4}{\rmdefault}{\mddefault}{\updefault}a}}}
\put(2176,-661){\makebox(0,0)[lb]{\smash{\SetFigFont{7}{8.4}{\rmdefault}{\mddefault}{\updefault}2}}}
\put(1726,-361){\makebox(0,0)[lb]{\smash{\SetFigFont{7}{8.4}{\rmdefault}{\mddefault}{\updefault}1}}}
\put(2401,-361){\makebox(0,0)[lb]{\smash{\SetFigFont{7}{8.4}{\rmdefault}{\mddefault}{\updefault}1}}}
\put(2251,-1036){\makebox(0,0)[lb]{\smash{\SetFigFont{7}{8.4}{\rmdefault}{\mddefault}{\updefault}a}}}
\end{picture}}
\parbox{0.4cm}{,}
\parbox{2cm}{\setlength{\unitlength}{0.00050000in}%
\begingroup\makeatletter\ifx\SetFigFont\undefined%
\gdef\SetFigFont#1#2#3#4#5{%
  \reset@font\fontsize{#1}{#2pt}%
  \fontfamily{#3}\fontseries{#4}\fontshape{#5}%
  \selectfont}%
\fi\endgroup%
\begin{picture}(762,864)(1726,-1063)
\thicklines
\put(2101,-811){\makebox(11.1111,16.6667){\SetFigFont{10}{12}{\rmdefault}{\mddefault}{\updefault}.}}
\put(2101,-811){\line( 0, 1){450}}
\multiput(2101,-361)(-11.84211,7.89474){20}{\makebox(11.1111,16.6667){\SetFigFont{7}{8.4}{\rmdefault}{\mddefault}{\updefault}.}}
\put(1801,-811){\line( 1, 0){675}}
\multiput(2101,-361)(11.84211,7.89474){20}{\makebox(11.1111,16.6667){\SetFigFont{7}{8.4}{\rmdefault}{\mddefault}{\updefault}.}}
\put(1876,-1036){\makebox(0,0)[lb]{\smash{\SetFigFont{7}{8.4}{\rmdefault}{\mddefault}{\updefault}a}}}
\put(2176,-1036){\makebox(0,0)[lb]{\smash{\SetFigFont{7}{8.4}{\rmdefault}{\mddefault}{\updefault}a+2}}}
\put(2176,-661){\makebox(0,0)[lb]{\smash{\SetFigFont{7}{8.4}{\rmdefault}{\mddefault}{\updefault}2}}}
\put(1726,-361){\makebox(0,0)[lb]{\smash{\SetFigFont{7}{8.4}{\rmdefault}{\mddefault}{\updefault}1}}}
\put(2401,-361){\makebox(0,0)[lb]{\smash{\SetFigFont{7}{8.4}{\rmdefault}{\mddefault}{\updefault}1}}}
\end{picture}}
\parbox{0.4cm}{,}
\end{center}

\noindent and the corresponding ones with $a+2$ replaced by $a-2$,
 (having the vertical strand on the $j$-th position)
 span all of $V(n,m)$. 
Indeed using the fusing matrices (which are invertible) we can
relate this system to the standard basis $L({\bf p})$. 

But now these are precisely the eigenvectors for $g_i$, because we
have the following relations:

\begin{center}
\parbox{0.5cm}{$\rho(g_i)$}
\parbox{2cm}{\setlength{\unitlength}{0.00050000in}%
\begingroup\makeatletter\ifx\SetFigFont\undefined%
\gdef\SetFigFont#1#2#3#4#5{%
  \reset@font\fontsize{#1}{#2pt}%
  \fontfamily{#3}\fontseries{#4}\fontshape{#5}%
  \selectfont}%
\fi\endgroup%
\begin{picture}(849,624)(964,-598)
\thicklines
\put(976,-586){\makebox(11.1111,16.6667){\SetFigFont{10}{12}{\rmdefault}{\mddefault}{\updefault}.}}
\put(976,-586){\line( 1, 0){825}}
\put(1351,-586){\line( 0, 1){375}}
\multiput(1351,-211)(-10.22727,10.22727){23}{\makebox(11.1111,16.6667){\SetFigFont{7}{8.4}{\rmdefault}{\mddefault}{\updefault}.}}
\multiput(1351,-211)(10.22727,10.22727){23}{\makebox(11.1111,16.6667){\SetFigFont{7}{8.4}{\rmdefault}{\mddefault}{\updefault}.}}
\put(1426,-511){\makebox(0,0)[lb]{\smash{\SetFigFont{7}{8.4}{\rmdefault}{\mddefault}{\updefault}k}}}
\put(1051,-211){\makebox(0,0)[lb]{\smash{\SetFigFont{7}{8.4}{\rmdefault}{\mddefault}{\updefault}1}}}
\put(1651,-211){\makebox(0,0)[lb]{\smash{\SetFigFont{7}{8.4}{\rmdefault}{\mddefault}{\updefault}1}}}
\end{picture}}
\parbox{0.5cm}{=}
\parbox{2cm}{\setlength{\unitlength}{0.00050000in}%
\begingroup\makeatletter\ifx\SetFigFont\undefined%
\gdef\SetFigFont#1#2#3#4#5{%
  \reset@font\fontsize{#1}{#2pt}%
  \fontfamily{#3}\fontseries{#4}\fontshape{#5}%
  \selectfont}%
\fi\endgroup%
\begin{picture}(849,924)(964,-598)
\thicklines
\put(976,-586){\makebox(11.1111,16.6667){\SetFigFont{10}{12}{\rmdefault}{\mddefault}{\updefault}.}}
\put(976,-586){\line( 1, 0){825}}
\put(1351,-586){\line( 0, 1){450}}
\multiput(1351,-136)(-12.50000,6.25000){13}{\makebox(11.1111,16.6667){\SetFigFont{7}{8.4}{\rmdefault}{\mddefault}{\updefault}.}}
\put(1201,-61){\line( 0, 1){150}}
\put(1201, 89){\line( 5, 3){375}}
\multiput(1351,-136)(13.23529,4.41176){18}{\makebox(11.1111,16.6667){\SetFigFont{7}{8.4}{\rmdefault}{\mddefault}{\updefault}.}}
\put(1576,-61){\line( 0, 1){150}}
\multiput(1576, 89)(-10.71429,10.71429){8}{\makebox(11.1111,16.6667){\SetFigFont{7}{8.4}{\rmdefault}{\mddefault}{\updefault}.}}
\multiput(1351,239)(-10.71429,10.71429){8}{\makebox(11.1111,16.6667){\SetFigFont{7}{8.4}{\rmdefault}{\mddefault}{\updefault}.}}
\put(1426,-511){\makebox(0,0)[lb]{\smash{\SetFigFont{7}{8.4}{\rmdefault}{\mddefault}{\updefault}k}}}
\put(1051,-211){\makebox(0,0)[lb]{\smash{\SetFigFont{7}{8.4}{\rmdefault}{\mddefault}{\updefault}1}}}
\put(1651,-211){\makebox(0,0)[lb]{\smash{\SetFigFont{7}{8.4}{\rmdefault}{\mddefault}{\updefault}1}}}
\end{picture}}
\parbox{0.5cm}{=}
\parbox{2cm}{\setlength{\unitlength}{0.00050000in}%
\begingroup\makeatletter\ifx\SetFigFont\undefined%
\gdef\SetFigFont#1#2#3#4#5{%
  \reset@font\fontsize{#1}{#2pt}%
  \fontfamily{#3}\fontseries{#4}\fontshape{#5}%
  \selectfont}%
\fi\endgroup%
\begin{picture}(849,870)(964,-598)
\thicklines
\put(1351,-211){\circle{300}}
\put(976,-586){\makebox(11.1111,16.6667){\SetFigFont{10}{12}{\rmdefault}{\mddefault}{\updefault}.}}
\put(976,-586){\line( 1, 0){825}}
\put(1351,-586){\line( 0, 1){225}}
\put(1351,-61){\line( 0, 1){150}}
\multiput(1351, 89)(-10.00000,10.00000){16}{\makebox(11.1111,16.6667){\SetFigFont{7}{8.4}{\rmdefault}{\mddefault}{\updefault}.}}
\multiput(1351, 89)(10.00000,10.00000){16}{\makebox(11.1111,16.6667){\SetFigFont{7}{8.4}{\rmdefault}{\mddefault}{\updefault}.}}
\put(1051,164){\makebox(0,0)[lb]{\smash{\SetFigFont{7}{8.4}{\rmdefault}{\mddefault}{\updefault}1}}}
\put(1651,164){\makebox(0,0)[lb]{\smash{\SetFigFont{7}{8.4}{\rmdefault}{\mddefault}{\updefault}1}}}
\put(1426,-511){\makebox(0,0)[lb]{\smash{\SetFigFont{7}{8.4}{\rmdefault}{\mddefault}{\updefault}k}}}
\put(1051,-286){\makebox(0,0)[lb]{\smash{\SetFigFont{7}{8.4}{\rmdefault}{\mddefault}{\updefault}1}}}
\put(1651,-286){\makebox(0,0)[lb]{\smash{\SetFigFont{7}{8.4}{\rmdefault}{\mddefault}{\updefault}1}}}
\end{picture}}
\parbox{2cm}{$= \; \delta(k;1,1) \; $}
\parbox{2cm}{\setlength{\unitlength}{0.00050000in}%
\begingroup\makeatletter\ifx\SetFigFont\undefined%
\gdef\SetFigFont#1#2#3#4#5{%
  \reset@font\fontsize{#1}{#2pt}%
  \fontfamily{#3}\fontseries{#4}\fontshape{#5}%
  \selectfont}%
\fi\endgroup%
\begin{picture}(849,624)(964,-598)
\thicklines
\put(976,-586){\makebox(11.1111,16.6667){\SetFigFont{10}{12}{\rmdefault}{\mddefault}{\updefault}.}}
\put(976,-586){\line( 1, 0){825}}
\put(1351,-586){\line( 0, 1){375}}
\multiput(1351,-211)(-10.22727,10.22727){23}{\makebox(11.1111,16.6667){\SetFigFont{7}{8.4}{\rmdefault}{\mddefault}{\updefault}.}}
\multiput(1351,-211)(10.22727,10.22727){23}{\makebox(11.1111,16.6667){\SetFigFont{7}{8.4}{\rmdefault}{\mddefault}{\updefault}.}}
\put(1426,-511){\makebox(0,0)[lb]{\smash{\SetFigFont{7}{8.4}{\rmdefault}{\mddefault}{\updefault}k}}}
\put(1051,-211){\makebox(0,0)[lb]{\smash{\SetFigFont{7}{8.4}{\rmdefault}{\mddefault}{\updefault}1}}}
\put(1651,-211){\makebox(0,0)[lb]{\smash{\SetFigFont{7}{8.4}{\rmdefault}{\mddefault}{\updefault}1}}}
\end{picture}}
\end{center}

\noindent This implies that the eigenvalues of $\rho(g_i)$ are $-A^{-3}$ and $A$, so that
shifting $\rho$ by a factor of $-A^{-1}$ will change them into $-1$ and
$A^{-4}$, as in the usual presentation of $H_n(q)$ with $q=A^{-4}$.

It can be proved (but this is beyond the scope of this note) that 
the representation $\tilde{\rho}$ is precisely the Hecke algebra 
representation $\pi_{\lambda}^{(r)}$ associated to the 
Young diagram $\lambda=[\frac{n+m}{2},\frac{n-m}{2}]$ as considered 
by Wenzl in \cite{We}. This is clear for $n=3, m=1$. The direct 
computational approach is somewhat cumbersome for the general case. 
However from \cite{Bl}  we can derive a still more general
equivalence between  the  
Hecke algebras representation associated to a Young diagram $\lambda$
with $N$ rows and that  arising 
in the previous construction for  the $SU(N)$-TQFT, where the label 
$m$ is replaced by the ``color'' $\lambda$. 
This follows from the irreducibility of the latter, using the 
technique from \cite{Bl}.

\end{rem}

\subsection{Proof of proposition \ref{pr3}}
In \cite{Jones2} a proof for  proposition \ref{pr3} is given for the case when 
$q=\exp(\frac{2\pi\sqrt{-1}}{r})$, but the argument generalizes easily
to all primitive roots of unity. We outline it below, for the sake of 
completeness.

It is known that $A_{\beta,3}$ is 
semi-simple and splits as 
$A_{\beta,3}= M_2({\bf C})\oplus {\bf C}$, for all $\beta\neq 1$ (see 
the theorem 2.8.5, p.98 from \cite{GHJ}). 
It suffices then to show that the images $\pi(g_1)$ and $\pi(g_2)$ 
in the factor $M_2({\bf C})$ generate an infinite group.  
Observe first that  $\pi(g_1)$ and $\pi(g_2)$
(and respectively $\tilde{\rho}(g_1)$ and $\tilde{\rho}(g_2)$)
do not commute with each other.  Thus 
the  $A_{\beta,3}$-module $V(3,1)$ is isomorphic to 
 the simple non-trivial factor $M_2({\bf C})$.
 As a consequence it suffices to see what happens with the
images of these two generators, when restricted to this summand. 
The $B_3$ representation on $M_2({\bf C})$  is also 
unitarizable when $q$ is a
root of unity according to proposition 3.2, p.257 from \cite{Jones2}.
Thus it makes sense to consider the images $\iota(g_1)$ and 
$\iota(g_2)$ in $SO(3)=U(2)/{\bf C}^*$. 
We have then the following decomposition in orthogonal projectors:
\[ \iota(g_i) = q e_i - (1-e_i), \]
so that the order of $\iota(g_i)$ in $SO(3)$ is $2r$ if $r$ is odd, 
$r/2$ if $r=2(4)$ and $r$ if $r=0(4)$, because $q$ is a primitive root
of unity of order $r$. As $r\neq 1$ these two elements cannot belong
to a cyclic or dihedral subgroup of $SO(3)$. But no other subgroups
have elements of order bigger than 5. Thus for $r=5,7,8,9$ or $r\geq
11$,  the image in $M_2({\bf C})$ of the subgroup generated by 
$g_1$ and $g_2$ is infinite. 

When $r=10$ we have to work within $B_4$, likewise to 
\cite{Jones2}, p.269. We already saw  
that the $A_{\beta,4}$-module $V(4,2)$ is 
irreducible and  identified therefore with the simple factor $M_3({\bf C})$ 
from $A_{\beta,4}$. 
Moreover the representation of $B_4$ 
on this factor was explicitly found out in \cite{Jones}, and it is 
the tensor product of the Burau and parity representations. 
It is  also shown that the  Burau representation
contains   elements of infinite order, for instance 
$g_1g_2g_3^{-1}$.   

We may wonder whether an element of infinite order in the image  
can be explicitly given for $r\neq 10$. Since we have to consider only the
matrices $\iota(B_n)$ in $SO(3)$, it is very likely that the 
element $g_1^{-1}g_2$ has infinite order. 

\newtheorem{lembb}{Remark}[section]
\begin{lembb}
Once we obtained the fact that the image of ${\cal M}_g$ is infinite
at a particular primitive root of unity, we may argue also as follows:
the Galois group $Gal(\overline{\bf Q};{\bf Q})$ acts on the set of
roots of unity, as well as on the entries of the matrices $\rho(x)$,
with $x\in {\cal M}_g$. It suffices to prove that the two actions
of $Gal(\overline{\bf Q};{\bf Q})$ are compatible to each other, in
order to conclude that the image group is infinite at all roots of
unity. This argument was pointed to me by Gregor Masbaum. 
\end{lembb}

\subsection{The RT version}
Lickorish \cite{L2} established the relationship among the invariants
obtained via the Temperley-Lieb algebra (basically those from
\cite{BHMV1}) $I(M,A)$ and the Reshetikhin-Turaev invariant $\tau_r(M)$
(see \cite{KiMe}), for closed oriented 3-manifolds $M$, as follows:
\[ I(M,
-\exp\frac{\pi\sqrt{-1}}{2r})=\exp(\frac{(6-3r)b_1(M)\pi\sqrt{-1}}{4r})
\tau_r(M), \]
where $b_1(M)$ is the first Betti number of the manifold $M$. 
Roughly speaking the two
invariants are the same up to a normalization factor. There are
however two associated TQFTs, still very close to each other:
\begin{enumerate}
\item The TQFT based on the Kauffman bracket, as described in
  \cite{BHMV2}, which arises in a somewhat canonical way; in fact any
  invariant of closed 3-manifolds extends to a TQFT via this procedure
(see \cite{BHMV2,F1} for details. The associated mapping class group
  representation we denote it by $\rho^K$. 
\item The TQFT based on the Jones polynomial, as described in
  \cite{KiMe}
(see also \cite{Ge}). The associated mapping class group
  representation we denote it by $\rho^J$, and may be computed using
  the definitions from conformal field theory like in \cite{MoSe}. 
A derivation of this representation, and the reconstruction of the
  invariant from it was first given by Kohno \cite{Kohno} (see also 
\cite{Tu1,Tu2,F5}). 
\end{enumerate}
The two representations are similar: the associated spaces on which
they act are naturally isomorphic. This means that in both theories
$W(F)$ has a distinguished basis given by labelings of 3-valent
graphs, with the same set of labels. Basically both theories are built
up using some variants of the quantum 6j-symbols:
\begin{enumerate}
\item in \cite{MV} these are identified with the tetrahedron
  coefficients, (see also \cite{KL}); the relationship with the usual
 6j-symbol (coming from representation theory) was outlined in
  \cite{Piu}. 
\item in the case of $\rho^J$ the 6j-symbols are coming from the
  representation theory of $U_q(sl_2)$ and where described in 
\cite{KiRe}. 
\end{enumerate}

Consider now the analogous subspace $V(n,m)=W(\Gamma'(n,m))$ of $W(F)$, as in
3.1. We have again an action of the braid group $B_n$ on $V$, but this
time the interpretation is no longer related to skein modules of the
ball. Here the graph $\Gamma$ is considered to be embedded in the
surface $F$, giving a rigid structure on $F$ \cite{F1,Wa}. This means
that there is a pants decomposition  $c$ of $F$ with the property that all
circles in  $c$ are transversal to $\Gamma$, the intersection of
$\Gamma$ with every trinion is the suspension of 3 points
(topologically, the space underlying the figure Y). Remark that $c$
and $\Gamma$ determine uniquely an identification of $F$ with a fixed 
and decomposed surface, up to an isotopy. 

This time twisting the strands of the labeled graphs in $L({\bf p})$ can
be expressed in terms of the data of conformal field theory (see
\cite{Kohno}). Specifically, we have:

\begin{center}
\parbox{2cm}{\setlength{\unitlength}{0.00050000in}%
\begingroup\makeatletter\ifx\SetFigFont\undefined%
\gdef\SetFigFont#1#2#3#4#5{%
  \reset@font\fontsize{#1}{#2pt}%
  \fontfamily{#3}\fontseries{#4}\fontshape{#5}%
  \selectfont}%
\fi\endgroup%
\begin{picture}(1074,789)(1564,-988)
\thicklines
\put(1876,-736){\line( 3, 4){387}}
\multiput(2251,-736)(-7.89474,11.84211){20}{\makebox(11.1111,16.6667){\SetFigFont{7}{8.4}{\rmdefault}{\mddefault}{\updefault}.}}
\multiput(2026,-436)(-7.89474,11.84211){20}{\makebox(11.1111,16.6667){\SetFigFont{7}{8.4}{\rmdefault}{\mddefault}{\updefault}.}}
\put(1576,-736){\line( 1, 0){1050}}
\put(1726,-361){\makebox(0,0)[lb]{\smash{\SetFigFont{7}{8.4}{\rmdefault}{\mddefault}{\updefault}a}}}
\put(2401,-361){\makebox(0,0)[lb]{\smash{\SetFigFont{7}{8.4}{\rmdefault}{\mddefault}{\updefault}b}}}
\put(1651,-961){\makebox(0,0)[lb]{\smash{\SetFigFont{7}{8.4}{\rmdefault}{\mddefault}{\updefault}c}}}
\put(2026,-961){\makebox(0,0)[lb]{\smash{\SetFigFont{7}{8.4}{\rmdefault}{\mddefault}{\updefault}d}}}
\put(2476,-961){\makebox(0,0)[lb]{\smash{\SetFigFont{7}{8.4}{\rmdefault}{\mddefault}{\updefault}e}}}
\end{picture}} 
\parbox{4cm}{$\; = \; 
\sum_{j} B_{dj}\left[\begin{array}{cc}
                     a & b \\
                     c & e 
                     \end{array}
                \right]\;$}
\parbox{2cm}{\setlength{\unitlength}{0.00050000in}%
\begingroup\makeatletter\ifx\SetFigFont\undefined%
\gdef\SetFigFont#1#2#3#4#5{%
  \reset@font\fontsize{#1}{#2pt}%
  \fontfamily{#3}\fontseries{#4}\fontshape{#5}%
  \selectfont}%
\fi\endgroup%
\begin{picture}(1074,744)(1564,-997)
\thicklines
\put(1576,-736){\line( 1, 0){1050}}
\put(1876,-736){\line( 0, 1){375}}
\put(2251,-736){\line( 0, 1){375}}
\put(1651,-361){\makebox(0,0)[lb]{\smash{\SetFigFont{7}{8.4}{\rmdefault}{\mddefault}{\updefault}a}}}
\put(2401,-361){\makebox(0,0)[lb]{\smash{\SetFigFont{7}{8.4}{\rmdefault}{\mddefault}{\updefault}b}}}
\put(1651,-961){\makebox(0,0)[lb]{\smash{\SetFigFont{7}{8.4}{\rmdefault}{\mddefault}{\updefault}c}}}
\put(2401,-961){\makebox(0,0)[lb]{\smash{\SetFigFont{7}{8.4}{\rmdefault}{\mddefault}{\updefault}e}}}
\put(2026,-961){\makebox(0,0)[lb]{\smash{\SetFigFont{7}{8.4}{\rmdefault}{\mddefault}{\updefault}j}}}
\end{picture}}
\end{center} 
where the matrix $B$  is the so-called braiding matrix. The braiding
matrix can be expressed in terms of the fusing matrix $F$ (see
\cite{Kohno,MoSe} by the following formula: 
\[ B_{ij}\left[\begin{array}{cc}
                     j_2 & j_3 \\
                     j_1 & j_4
                     \end{array}
                \right]
= (-1)^{j_1+j_4-i-j}/2 \exp(\pi\sqrt{-1}(\Delta_{j_1}+\Delta_{j_4}-
\Delta_{i}-\Delta_{j})) 
F_{ij}\left[\begin{array}{cc}
                     j_1 & j_3 \\
                     j_2 & j_4
                     \end{array}
                \right], \]
where 
\[ \Delta_{j} =\frac{j(j+1)}{4r}. \]
We use the same set of labels for the graphs, namely integers running 
 from $1$ to
$2r-2$ as before,  instead of the traditional  half-integer labels from
\cite{KiRe,Kohno,KL}. Set also $q=\exp\frac{2\pi\sqrt{-1}}{r}$, and 
$[n]=\frac{q^n-q^{-n}}{q-q^{-1}}$.

The natural choice for the fusing matrix $F$ is (see
\cite{Kohno},p.213-214,\cite{TuVi}):
\[ F_{ij}\left[\begin{array}{cc}
                     j_1 & j_3 \\
                     j_2 & j_4
                     \end{array}
                \right] =
\left\{ \begin{array}{ccc}
                     j_1 & j_3 & i \\
                     j_2 & j_4 & j
                     \end{array}
                \right \}_{KR}, \]
where $\{,\}_{KR}$ denotes the quantum 6j-symbols of Kirillov and
Reshetikhin. 

Using the computations from \cite{KiRe}, and those from 2.3 we find
that the only non-trivial braiding matrix for $j_2=j_3=1$ is that 
with $j_1=j_4$, and its value is therefore:
\[ B\left[\begin{array}{cc}
                     1 & 1 \\
                     a & a
                     \end{array}
                \right]= 
 \left(\begin{array}{cc}
                    -q^{a+\frac{1}{4}}\left(\frac{[a]}{[2][a+1]}\right)^{1/2} 
& -q^{-\frac{1}{4}}\left(\frac{[a+2]}{[2][a+1]}\right)^{1/2} \\
   -q^{-\frac{1}{4}}\left(\frac{[a+2]}{[2][a+1]}\right)^{1/2}
     & -q^{-a-\frac{3}{4}}\left(\frac{[a]}{[2][a+1]}\right)^{1/2}
                     \end{array}
                \right). \]
Notice that the braiding matrices arising in conformal field theory 
 were previously computed by Tsuchyia
and Kanie in \cite{TK}. Their result, used however a different
normalization and the matrices are not identical, but equivalent up to
a power of $q$. In fact, in our case, the representation
$q^{1/4}\rho^J$ is also equivalent to $\tilde{\rho}$, for 
$n=3,4$. 
As an immediate consequence the representation $\rho^J$ has an
infinite image, too, under the same condition as $\rho^K$. 
This ends the proof of the main theorem.

\bibliographystyle{plain}

\end{document}